\documentclass{amsart}

\usepackage[T1]{fontenc}
\usepackage{times}
\usepackage{amssymb,verbatim}
\theoremstyle{plain}

\usepackage[utf8]{inputenc}
\usepackage[english]{babel}
\usepackage{amssymb}
\usepackage{amsthm}
\usepackage{graphics}
\usepackage{amsmath}
\usepackage{amstext}
\usepackage{multirow}
\usepackage{graphicx}
\usepackage{amsmath}
\usepackage{amssymb}
\usepackage{amsthm}
\usepackage{amstext}
\usepackage{mathrsfs}
\usepackage{amscd}
\usepackage{mathtools} 
\usepackage[
    backref=page
]{hyperref} 
\usepackage{tikz}
\usetikzlibrary{matrix}
\usepackage[all]{xy}
\usepackage{algorithm}
\usepackage{algorithmic}
\usepackage{comment}

\newtheorem{thm}{Theorem}[section]
\newtheorem*{conj*}{Conjecture}

\newtheorem{lemma}[thm]{Lemma}
\newtheorem{prop}[thm]{Proposition}
\newtheorem{cor}[thm]{Corollary}

\newtheorem{THM}{Theorem}

\theoremstyle{remark}

\newtheorem{example}[thm]{Example}
\newtheorem{remark}[thm]{Remark}

\newcommand{\mb}{\mathbb}
\newcommand{\mc}{\mathcal}

\newcommand{\C}{\mb C}
\newcommand{\Pj}{\mb P}

\newcommand{\Q}{\mb Q}

\newcommand{\Proj}{\mathbb{P}} 

\newcommand{\F}{\mc F}

\DeclareMathOperator{\sing}{Sing}

\DeclareMathOperator{\Pic}{Pic}
\DeclareMathOperator{\GSV}{GSV}
\DeclareMathOperator{\QGSV}{GSV^{\Q}}
\DeclareMathOperator{\CS}{CS}
\DeclareMathOperator{\QCS}{CS^{\Q}}
\DeclareMathOperator{\Cl}{Cl}

\newcommand{\NF}{{N_\F}}

\numberwithin{equation}{section}

\usepackage{mathtools}
\usepackage{ifmtarg}

\addtocounter{section}{0}             
\numberwithin{equation}{section}       

\title{Holomorphic Foliations with a Single Singularity on Orbifold Projective Surfaces}

\author[R. Lizarbe]{Ruben Lizarbe}
\address{UERJ, Universidade do Estado do Rio de Janeiro, Rua S\~ao Francisco Xavier, 524, Maracan\~a, 20550-900, Rio de Janeiro, Brazil. }
\email{ruben.monje@ime.uerj.br}

\author[L. Puchuri]{Liliana Puchuri}
\address{Departamento de Ciencias, Secci\'on Matem\'aticas, Pontificia Universidad Cat\'olica del Per\'u, Av. Universitaria 1801, San Miguel, Lima, Per\'u}
\email{lpuchuri@pucp.pe}

\author[R. Rosas]{Rudy Rosas}
\address{Departamento de Ciencias, Secci\'on Matem\'aticas, Pontificia Universidad Cat\'olica del Per\'u, Av. Universitaria 1801, San Miguel, Lima, Per\'u}
\email{rudy.rosas@pucp.edu.pe}


\begin{document}

\begin{abstract}
We investigate holomorphic foliations on orbifold projective surfaces with cyclic divisor class group, assuming that the foliation has a unique singular point of algebraic multiplicity one. If the normal degree of the foliation is sufficiently large, then we prove that it admits no invariant algebraic curve. Furthermore, we present several examples showing that the bound obtained is optimal.

\end{abstract}

\maketitle

\setcounter{tocdepth}{1}

\tableofcontents

\section{Introduction}

A holomorphic foliation on a projective surface need not leave any algebraic
curve invariant, and deciding whether a given foliation does is a classical
and delicate problem. The subject begins with Jouanolou \cite{Jou}, who
proved that a generic foliation of degree $d\ge2$ on the complex projective
plane $\Proj^2$ admits no invariant algebraic curve, and who exhibited, for
each $d\ge2$, the explicit foliation induced by the homogeneous polynomial
$1$-form
\[
\Omega_d=(x^dz-y^{d+1})\,dx+(xy^d-z^{d+1})\,dy+(zy^d-x^{d+1})\,dz.
\]

Genericity results of this type have since been extended in several
directions: a different proof for $\Proj^2$ was given in \cite{Jorge}, and
\cite{Cou} shows that a generic foliation of sufficiently large normal
degree on an arbitrary smooth projective surface admits no invariant
algebraic curve; analogous statements hold on weighted projective planes
\cite{lizarbe2023density}. All of these are existence theorems: they show
that the foliations in question form a dense family, but they provide
neither a criterion for testing a given foliation nor a recipe for
producing new examples. Jouanolou's $\Omega_d$ remains the model of what
one would like to have instead.

One way around this obstacle is to single out foliations with a unique
singular point --- a hypothesis rigid enough to be verified in practice,
yet strong enough to force the absence of invariant curves; see
\cite{alcantara2018foliations, alcantara2020foliations,
cerveau2010geometrie, FernandezPuchuriRosas2021}. Our starting point is a
theorem of Fernández, Puchuri, and Rosas
\cite{FernandezPuchuriRosas2021}: a foliation on $\Proj^2$ whose unique
singularity has algebraic multiplicity one admits no invariant algebraic
curve. Both hypotheses can be checked by direct computation, so the theorem
also serves to construct explicit examples.

Our aim is to carry this result beyond the smooth setting, to orbifold
projective surfaces --- normal projective surfaces whose singularities are
quotient singularities. Weighted projective planes are the motivating case,
and already there the classical arguments break down: the ambient surface
is itself singular, so the Milnor number and the Baum--Bott,
Gómez-Mont--Seade--Verjovsky and Camacho--Sad indices must be replaced by
orbifold counterparts, defined through local uniformizations and taking
rational rather than integer values. Section~\ref{F:Formulas} constructs
these substitutes and proves the index formulas they satisfy; they are the
technical backbone of everything that follows.

Our first result assumes only that the divisor class group is cyclic.

\begin{THM}\label{THM:A}
Let $X$ be an orbifold projective surface with $\Cl(X)\cong\mathbb{Z}$ and
$c_2(X)>1$. Let $H$ be the ample generator of $\Cl(X)$, and let $\F$ be a
foliation on $X$ with a unique singularity at $x\in X$. Let
$d=\deg(N_{\F})$, so that $N_{\F}\equiv dH$, and assume that the algebraic
multiplicity of $\F$ at $x$ is one. Then there exists a constant $C$,
depending only on $X$, such that if $d>C$, then $\F$ admits no invariant
algebraic curve.
\end{THM}

An explicit value of $C$ is given in Theorem~\ref{T:classone}. When $X$ is
smooth or toric, we further show that the strict inequality $d>C$ can be
slightly relaxed, and that the bound is optimal in several examples; the
following theorem confirms this for a family of weighted projective planes.

\begin{THM}\label{THM:B}
Let $\F$ be a foliation on $\Pj(1,1,k)$, $k>1$, of normal degree $d$ with a
unique singular point $p$. Assume that $\F$ has algebraic multiplicity one
at $p$, that is, $m_p(\F)=1$. Then the following assertions hold:
\begin{enumerate}
\item If $d\ge 2k+2$, then $\F$ admits no invariant algebraic curve.
\item If $d=2k+1$, then $\F$ admits no invariant algebraic curve when $p$
is a saddle-node singularity, whereas it may admit one when $p$ is
nilpotent.
\end{enumerate}
\end{THM}

The second alternative marks a genuine departure from the smooth case. On
$\Proj^2$, a foliation with a unique nilpotent singularity never admits an
invariant algebraic curve \cite[Theorem 1]{cerveau2010geometrie},
\cite[Theorem 1]{FernandezPuchuriRosas2021}; on $\Pj(1,1,k)$ with $k>1$
this fails, and Remark~\ref{R:nilpontentfirstintegral} exhibits a
counterexample.

The bound of Theorem~\ref{THM:B} is sharp, and the two theorems below
explain why: they classify the foliations that sit exactly at the
threshold. We begin with normal degree $2k+1$.

\begin{THM}\label{THM:C}
Let $\F$ be a foliation on $\Pj(1,1,k)$ of normal degree $2k+1$ having a
unique singular point $p$. Then exactly one of the following alternatives
occurs:
\begin{enumerate}
\item $p=[0:0:1]$, $m_p(\F)=k$, and $\F$ is a Riccati foliation. In
particular, $\F$ admits an invariant algebraic curve.

\item After an automorphism of $\Pj(1,1,k)$, we may assume that
$p=[1:0:0]$. In suitable affine coordinates $(y,z)$, the foliation $\F$ is
induced by one of the following $1$-forms:
\begin{enumerate}
\item
\[
\omega=
\Bigl(
y^{2k-1}
+z\sum_{j=2}^{k}a_jy^{j-1}
+kz^2
\Bigr)\,dy
-yz\,dz,
\qquad a_j\in\C .
\]

\item
\[
\omega=
\Bigl(
(z-y^{k+1})
\sum_{j=1}^{k-1}b_jy^{j-1}
+kzy^k
\Bigr)\,dy
+(z-y^{k+1})\,dz,
\qquad b_j\in\C .
\]

\item
\[
\omega=
\Bigl(
(z-y^{k+1}-zy)
\sum_{j=1}^{k-1}b_jy^{j-1}
+kz(y^k+z)
\Bigr)\,dy
+(z-y^{k+1}-zy)\,dz,
\qquad b_j\in\C .
\]
\end{enumerate}
\end{enumerate}
\end{THM}

The two alternatives are distinguished by the position of the singularity:
it lies at $[0:0:1]$ exactly in the Riccati case, and at $[1:0:0]$
otherwise. In the latter case, form {\rm (a)} yields a singularity of
algebraic multiplicity $2$, whereas forms {\rm (b)} and {\rm (c)} yield a
saddle-node singularity when $b_1\neq 0$ and a nilpotent one when $b_1=0$.
Together with Theorem~\ref{THM:B}, this accounts for every foliation of
normal degree $2k+1$ with a single singular point.

The remaining case, normal degree $2k+2$, admits a shorter description.

\begin{THM}\label{THM:D}
Let $\F$ be a foliation on $\Pj(1,1,k)$ of normal degree $2k+2$, with a
unique singularity at $p=[0:0:1]$. Then exactly one of the following holds:
\begin{enumerate}
\item $m_p(\F)=k+1$, and $\F$ is a Riccati foliation. In particular, $\F$
admits at least one invariant algebraic curve; or

\item $m_p(\F)=1$, and $p$ is a saddle-node singularity. In this case, $\F$
admits no invariant algebraic curve.
\end{enumerate}
\end{THM}

\subsection{Organization of the paper}

Section~\ref{S:preliminaries} recalls the construction of weighted
projective planes and of the foliations they carry, together with the
notions of normal degree and algebraic multiplicity.
Section~\ref{F:Formulas} is the technical core of the paper: there we
introduce the orbifold Milnor number and the orbifold Baum--Bott,
Gómez-Mont--Seade--Verjovsky and Camacho--Sad indices, we establish the
formulas relating them, and we derive the degree bound of
Proposition~\ref{P:bounddegree}. Section~\ref{S:classone} contains the
proof of Theorem~\ref{THM:A}, in the explicit form of
Theorem~\ref{T:classone}; the saddle-node and the nilpotent case are
treated separately, and the smooth and toric surfaces are singled out,
since for them the bound can be slightly improved.
Section~\ref{S:applications} specializes the results to weighted
projective planes and contains the proofs of Theorems~\ref{THM:B},
\ref{THM:C} and~\ref{THM:D}, together with explicit bounds for
$\Pj(1,b,c)$ and $\Pj(2,3,5)$.

\subsection{Acknowledgments}

\subsection*{Use of Large Language Models}

During the development and preparation of this article, the authors used (OpenAI; GPT-5.6 Luna) and Claude (Anthropic; Claude Fable 5) 
as interactive tools to explore proof strategies and to assist with 
drafting and revision. The authors take full responsibility for the 
content of the article.

\section{Orbifold projective surfaces and Holomorphic foliations}
\label{S:preliminaries}

This section recalls the basic notions concerning orbifold projective
surfaces and holomorphic foliations that will be used throughout the
paper. We refer to \cite{Mumford1983} for the orbifold viewpoint,
to \cite{Dolgachev} for weighted projective spaces, and to
\cite{Brun2, Lizarbe14, lizarbe2023density} for foliations in this
setting. We also use the description of foliations on normal
varieties with quotient singularities developed in
\cite{Druel-Ou2020}.

\subsection{Orbifold Projective Surfaces}

Let $X$ be a normal projective surface with quotient singularities.
Following Mumford's approach to orbifold structures
\cite[Section~2]{Mumford1983}, we study $X$ locally through finite
quotient charts
\[
p_\alpha:\widetilde U_\alpha\longrightarrow U_\alpha\subset X,
\]
where $\widetilde U_\alpha$ is smooth and $U_\alpha$ is obtained
as a quotient by a finite group. Geometric objects on $X$ can
therefore be studied by lifting them to these smooth charts and
requiring compatibility with the corresponding group actions.

Since $X$ may be singular, its cotangent and tangent sheaves need
not be locally free at every point. We use the reflexive cotangent
sheaf
\[
\Omega_X^{[1]}:=(\Omega_X^1)^{**}
\]
and the reflexive tangent sheaf
\[
T_X:=\mathcal Hom_{\mathcal O_X}
\bigl(\Omega_X^{[1]},\mathcal O_X\bigr),
\]
where $(-)^{**}$ denotes the double dual.

A rank-one reflexive sheaf $\mathscr L$ on $X$ is called a
$\mathbb Q$-line bundle if there exists an integer $m>0$ such that
its reflexive power
\[
\mathscr L^{[m]}
:=
\bigl(\mathscr L^{\otimes m}\bigr)^{**}
\]
is an invertible sheaf. Since quotient singularities are
$\mathbb Q$-factorial, Weil divisors on $X$ are $\mathbb Q$-Cartier.
In particular, the canonical sheaf
\[
\omega_X:=\bigl(\Omega_X^2\bigr)^{**}
\]
determines a $\mathbb Q$-Cartier canonical divisor $K_X$.

\subsection{Holomorphic Foliations}

A holomorphic foliation $\mathcal F$ of rank one on $X$ is given
by a coherent subsheaf
\[
T_{\mathcal F}\subset T_X
\]
of rank one such that $T_X/T_{\mathcal F}$ is torsion-free.
Equivalently, $T_{\mathcal F}$ is saturated in $T_X$.

On a smooth orbifold chart, the lifted foliation can be described
locally by a holomorphic vector field $v$ generating its tangent
direction. Dually, it can be described by a holomorphic $1$-form
$\omega$ annihilating that direction, with values in the
appropriate line bundle. Thus, locally,
\[
\omega(v)=0.
\]
These descriptions are compatible with the finite group actions
defining the quotient charts.

The canonical sheaf of the foliation is determined by
\[
\mathcal O_X(-K_{\mathcal F})
\simeq
\bigl(\det T_{\mathcal F}\bigr)^{**}.
\]
The normal sheaf of the foliation is the rank-one reflexive sheaf
\[
N_{\mathcal F}
:=
\bigl(T_X/T_{\mathcal F}\bigr)^{**}.
\]
At the level of divisor classes, the determinant relation gives
\[
K_{\mathcal F}
\sim_{\mathbb Q}
K_X+c_1(N_{\mathcal F}).
\]
These sheaves and divisor classes are understood in the reflexive
and $\mathbb Q$-Cartier sense appropriate to quotient singularities.

\subsection{Invariant Curves on Orbifold Surfaces}

Let $X$ be a normal projective surface with quotient singularities,
and let $\mathcal F$ be a holomorphic foliation of rank one on $X$.
An irreducible algebraic curve $C\subset X$ is said to be
\emph{invariant by $\mathcal F$} if, at a general smooth point
$p\in C\cap X_{\mathrm{reg}}$ outside the singular set of the
foliation, the tangent line $T_pC$ coincides with the line
$T_{\mathcal F,p}$ defined by the foliation.

This condition can be expressed locally on a smooth orbifold chart
\[
p:\widetilde U\longrightarrow U\subset X.
\]
Let $\widetilde{\mathcal F}$ and $\widetilde C$ denote the lifts
of the foliation and the curve, respectively. Suppose that
$\widetilde C$ is locally defined by a holomorphic equation
$f=0$. The lifted foliation can be described either by a
holomorphic vector field $v$ generating its tangent direction or
by a holomorphic $1$-form $\omega$ annihilating that direction,
so that
\[
\omega(v)=0.
\]

The curve $\widetilde C$ is invariant if and only if the vector
field is tangent to it. In terms of its local equation, this is
equivalent to
\[
v(f)\in(f),
\]
where $(f)$ denotes the ideal generated by $f$ in
$\mathcal O_{\widetilde U}$. Equivalently, using the $1$-form
description, the invariance condition is
\[
\omega\wedge df\in
(f)\,\Omega_{\widetilde U}^{2}.
\]
These conditions express the tangency of the lifted foliation
to the lifted curve. Since the local descriptions are compatible
with the finite group actions defining the orbifold charts, the
invariance condition is independent of the choice of chart.

We now introduce weighted projective planes, which provide an important class of examples of orbifold projective surfaces and play a central role in this work. We recall their basic properties and describe the holomorphic foliations they carry.

\subsection{Weighted Projective Planes}

Throughout this subsection, let $a,b,c\in\mathbb N$ be pairwise
coprime integers with $1\leq a\leq b\leq c$. The weighted projective
plane is
\[
\mathbb P(a,b,c)
=
(\mathbb C^3\setminus\{0\})/\mathbb C^*,
\]
where $\mathbb C^*$ acts by
\[
\lambda\cdot(x,y,z)
=
(\lambda^a x,\lambda^b y,\lambda^c z),
\qquad \lambda\in\mathbb C^*.
\]
We write $[x:y:z]$ for the class of $(x,y,z)$. When $a=b=c=1$,
we recover the usual projective plane
$\mathbb P^2=\mathbb P(1,1,1)$. In general,
$\mathbb P(a,b,c)$ is a normal projective surface with quotient
singularities, so the preceding definitions apply.

The standard charts
\[
U_x=\{x\neq0\},\qquad
U_y=\{y\neq0\},\qquad
U_z=\{z\neq0\}
\]
cover $\mathbb P(a,b,c)$. For example,
\[
U_z\simeq\mathbb C^2/\mu_c,
\]
where $\mu_c$ is the group of $c$-th roots of unity acting by
\[
(x,y)\longmapsto(\xi^a x,\xi^b y).
\]
Consequently, the singularities can occur only at the coordinate
points
\[
[1:0:0],\qquad [0:1:0],\qquad [0:0:1],
\]
although some of these points may be smooth.

The finite morphism
\[
\varphi:\mathbb P^2\longrightarrow\mathbb P(a,b,c),
\qquad
[X:Y:Z]\longmapsto[X^a:Y^b:Z^c],
\]
realizes $\mathbb P(a,b,c)$ as a finite quotient of $\mathbb P^2$
by the coordinatewise action of
$\mu_a\times\mu_b\times\mu_c$, with the projective action
understood modulo its kernel.

To simplify notation, we write
\[
\Proj:=\mathbb P(a,b,c).
\]

\subsection{Holomorphic Foliations on Weighted Projective Planes}

An algebraic curve in $\Proj$ is defined by a nonconstant
quasi-homogeneous polynomial $F\in\mathbb C[x,y,z]$:
\[
C=V(F)
=
\bigl\{
[x:y:z]\in\Proj\ ;\ F(x,y,z)=0
\bigr\}.
\]
The polynomial $F$ is quasi-homogeneous of degree $m$ if
\[
F(\lambda^a x,\lambda^b y,\lambda^c z)
=
\lambda^m F(x,y,z)
\qquad\text{for every }\lambda\in\mathbb C^*.
\]

In the convention adopted here, a holomorphic foliation
$\mathcal F$ of normal degree $d$ on $\Proj$ is represented on
the affine cone by a quasi-homogeneous $1$-form
\[
\Omega=A\,dx+B\,dy+C\,dz,
\]
where $A$, $B$ and $C$ are quasi-homogeneous polynomials of
weighted degrees $d-a$, $d-b$ and $d-c$, respectively. The form
has weighted degree $d$ and satisfies
\[
i_R\Omega=0,
\]
where
\[
R
=
ax\frac{\partial}{\partial x}
+
by\frac{\partial}{\partial y}
+
cz\frac{\partial}{\partial z}
\]
is the weighted radial vector field. Equivalently,
\[
axA+byB+czC=0.
\]
This condition ensures that the form is horizontal with respect
to the weighted $\mathbb C^*$-action. Together with its
quasi-homogeneity, it gives the corresponding foliation on
$\Proj$, interpreted in terms of reflexive sheaves at the quotient
singularities.

The singular set of $\mathcal F$ is
\[
\operatorname{Sing}(\mathcal F)
=
\bigl\{
[x:y:z]\in\Proj\ ;\
A(x,y,z)=B(x,y,z)=C(x,y,z)=0
\bigr\}.
\]
We assume that this set is discrete.

Let $H$ denote the $\mathbb Q$-Cartier divisor class corresponding
to $\mathcal O_{\Proj}(1)$. The weighted Euler sequence gives
\[
K_{\Proj}\sim -(a+b+c)H.
\]
Moreover, for a foliation of normal degree $d$ defined by $\Omega$,
the normal sheaf is
\[
N_{\mathcal F}\simeq\mathcal O_{\Proj}(d),
\]
where this notation is understood in the reflexive sense. Hence
\[
K_{\mathcal F}
\sim_{\mathbb Q}
K_{\Proj}+c_1(N_{\mathcal F})
\sim (d-a-b-c)H.
\]
In particular,
\[
\bigl(\det T_{\mathcal F}\bigr)^{**}
\simeq\mathcal O_{\Proj}(a+b+c-d).
\]
Since $T_{\mathcal F}$ has rank one, this yields
\[
T_{\mathcal F}
\simeq\mathcal O_{\Proj}(a+b+c-d)
\]
as a rank-one reflexive sheaf. These formulas are understood in
the $\mathbb Q$-Cartier and reflexive sense, since
$\mathcal O_{\Proj}(1)$ need not be invertible at the quotient
singularities.

\subsection{Invariant Curves on Weighted Projective Planes}

Let $\mathcal F$ be a foliation on $\Proj$ represented by
\[
\Omega=A\,dx+B\,dy+C\,dz,
\]
and let $C=V(F)$ be an irreducible algebraic curve defined by
a quasi-homogeneous polynomial $F$. By the definition of an
invariant curve on an orbifold surface, $C$ is invariant by
$\mathcal F$ if the lifted curve is tangent to the lifted
foliation at general points of a smooth local orbifold chart.

On the affine cone, this tangency condition is expressed by
\[
\Omega\wedge dF=0
\qquad\text{along }V(F).
\]
Equivalently, the coefficients of $\Omega\wedge dF$ are divisible
by $F$, so that
\[
F\mid\Omega\wedge dF.
\]
This gives an algebraic criterion for invariance in terms of
the quasi-homogeneous equation of the curve and the $1$-form
defining the foliation.

For example, consider
\[
X=\mathbb P(1,1,k),\qquad k\geq1,
\]
with coordinates $[x:y:z]$ of weights $1,1,k$. The $1$-form
\[
\Omega=y\,dx-x\,dy
\]
has weighted degree $2$ and satisfies
\[
i_R\Omega=xy-yx=0,
\qquad
R=x\frac{\partial}{\partial x}
+y\frac{\partial}{\partial y}
+kz\frac{\partial}{\partial z}.
\]
It therefore induces a foliation on $X$ in the reflexive-sheaf
sense. With the degree convention above, this foliation has
normal degree $2$.

For each $[u:v]\in\mathbb P^1$, consider the curve
\[
C_{[u:v]}=V(ux-vy).
\]
Writing $F=ux-vy$, we compute
\[
\begin{aligned}
\Omega\wedge dF
&=(y\,dx-x\,dy)\wedge(u\,dx-v\,dy)\\
&=(ux-vy)\,dx\wedge dy\\
&=F\,dx\wedge dy.
\end{aligned}
\]
Hence $F\mid\Omega\wedge dF$, and every curve $C_{[u:v]}$
is invariant by the foliation. In particular, the coordinate
curves $x=0$ and $y=0$ are invariant.

This example illustrates how invariant curves on weighted
projective planes can be studied using quasi-homogeneous
equations and the differential forms defining foliations.

\section{Orbifold Indices and Formulas}\label{F:Formulas}

Every classical local index of a foliation --- the Milnor number and the
Baum--Bott, Gómez-Mont--Seade--Verjovsky and Camacho--Sad indices --- has an
orbifold counterpart, and all four are obtained by the same recipe: lift the
foliation to a local uniformization and average the classical index over the
fibre. The resulting quantities are rational rather than integral, and they
satisfy exact analogues of the classical index formulas. Establishing this is
the purpose of the present section; the arguments of
Sections~\ref{S:classone} and~\ref{S:applications} rest entirely on it.

Throughout, $X$ denotes a normal projective surface with quotient
singularities and $\F$ a foliation on $X$. We fix the following notation once
and for all. Given a point $x\in X$, there exist an analytic neighbourhood $U$
of $x$, a (not necessarily connected) smooth complex surface $V$, and a finite
Galois holomorphic map
\[
p:V\longrightarrow U
\]
which is étale outside the singular locus of $U$; we call $\pi$ a \emph{local
uniformization} at $x$. We denote by $\widetilde{\F}$ the foliation induced by
$\F|_U$ on $V$ and, whenever a curve $C\subset X$ is given, by $\widetilde{C}$
the preimage $p^{-1}(C\cap U)$. None of the definitions below depends on the
choice of $p$.

The \emph{orbifold Milnor number} of $\F$ at $x$ is
\[
\mu^{\Q}(\F,x)
=
\frac{1}{\deg(p)}
\sum_{y\in p^{-1}(x)}
\mu(\widetilde{\F},y),
\]
where $\mu(\widetilde{\F},y)$ is the usual Milnor number of $\widetilde{\F}$ at
$y$, and the \emph{total orbifold Milnor number} of $\F$ is
\[
\mu^{\Q}(\F)
=
\sum_{x\in\sing(\F)}
\mu^{\Q}(\F,x).
\]

\begin{prop}\label{P:Q-milnor}
Let $\F$ be a foliation on a normal projective surface $X$ with quotient
singularities. Then
\[
\mu^{\Q}(\F)
=
c_2(X)-c_1(T_X)\cdot N_{\F}+N_{\F}^2
=
c_2(X)-T_{\F}\cdot N_{\F}.
\]
\end{prop}

\begin{proof}
By the Bott-type residue formula for holomorphic sections of orbifold
vector bundles \cite[Theorem~2.1, p.~2895]{CorreaRodriguezPenaSoares},
we have
\[
\mu^{\mathbb Q}(\mathcal F)
=c_2(T_X\otimes K_{\mathcal F}).
\]
Since $K_{\mathcal F}$ and $N_{\mathcal F}$ are $\mathbb Q$-line
bundles, the Chern class formula gives
\[
\begin{aligned}
\mu^{\mathbb Q}(\mathcal F)
&=c_2(X)-c_1(T_X)\cdot N_{\mathcal F}
  +N_{\mathcal F}^2\\
&=c_2(X)-T_{\mathcal F}\cdot N_{\mathcal F},
\end{aligned}
\]
where we use
$c_1(K_{\mathcal F})=-c_1(T_X)+N_{\mathcal F}$
and $c_1(T_{\mathcal F})=c_1(T_X)-N_{\mathcal F}$.
This proves the proposition.
\end{proof}


The \emph{orbifold Baum--Bott index} of $\F$ at $x$ is defined in the same way,
\[
BB^{\Q}(\F,x)
=
\frac{1}{\deg(p)}
\sum_{y\in p^{-1}(x)}
BB(\widetilde{\F},y),
\]
where $BB(\widetilde{\F},y)$ is the Baum--Bott index of $\widetilde{\F}$ at $y$.

\begin{prop}\label{P:Q-BB}
Let $\F$ be a foliation on $X$. Then
\[
c_1(\NF)^2
=
\sum_{x\in X}
BB^{\Q}(\F,x).
\]
\end{prop}

\begin{proof}
The statement follows from \cite[Proposition~3.9, p.~15]{Druel-Ou2020}.
\end{proof}

Let now $C\subset X$ be a complete $\F$-invariant curve. The \emph{orbifold
Gómez-Mont--Seade--Verjovsky index} of $\F$ along $C$ at $x$ is
\[
\QGSV(\F,C,x)
=
\frac{1}{\deg(p)}
\sum_{y\in p^{-1}(x) \cap \widetilde{C}}
\GSV(\widetilde{\F},\widetilde{C},y),
\]
and, summing over the singular points,
\[
\QGSV(\F,C)
=
\sum_{x\in\sing(\F) \cap C}
\QGSV(\F,C,x).
\]
Likewise, the \emph{orbifold Camacho--Sad index} of $\F$ along $C$ at $x$ is
\[
\QCS(\F,C,x)
=
\frac{1}{\deg(p)}
\sum_{y\in p^{-1}(x) \cap \widetilde{C}}
\CS(\widetilde{\F},\widetilde{C},y).
\]
Here $\GSV(\widetilde{\F},\widetilde{C},y)$ and $\CS(\widetilde{\F},\widetilde{C},y)$
denote the classical Gómez-Mont--Seade--Verjovsky and Camacho--Sad indices of
$\widetilde{\F}$ along $\widetilde{C}$ at $y$.

\begin{prop}\label{P:Q-GVS-CS}
Let $\F$ be a foliation on $X$, and let $C\subset X$ be a complete $\F$-invariant curve. Then
\[
\NF\cdot C-C^{2}
=
\QGSV(\F,C)
\]
and
\[
C^{2}
=
\sum_{x\in C}
\QCS(\F,C,x).
\]
\end{prop}

\begin{proof}
The first identity follows from the proof of
\cite[Proposition~3.12, p.~18]{Druel-Ou2020}, while the second is precisely
\cite[Proposition~3.12, p.~18]{Druel-Ou2020}.
\end{proof}

%
%
%
%

One further index will be needed. Let \(X\) be a normal projective
surface with quotient singularities, and let \(\F\) be a foliation
on \(X\). Let \(B\) be a local branch of an \(\F\)-invariant
curve \(C\) at a point \(x\).

If \(x\in\operatorname{Sing}(X)\), let \(G_x\) be the local
uniformizing group and let
\[
\pi:(\mathbb C^2,0)\longrightarrow (X,x)
\]
be the corresponding quotient map. Denote by
\(\widetilde{\F}\) the lifted foliation. We define the
\emph{orbifold relative Milnor number} of \(\F\) along \(B\) by
\[
\mu^{\mathbb Q}(\F,B,x)
=
\frac{\mu(\widetilde{\F},\pi^{-1}(B),0)}{|G_x|}.
\]
If \(x\) is a smooth point of \(X\), we set
\[
\mu^{\mathbb Q}(\F,B,x)=\mu(\F,B,x).
\]

\begin{prop}\label{P:Q-milnorrelativo}
Let \(X\) be a normal projective surface with quotient
singularities, let \(\F\) be a foliation on \(X\), and let
\(C\subset X\) be a reduced \(\F\)-invariant curve. 
Then
\[
\sum_{x\in C\cap\operatorname{Sing}(\F)}
\ \sum_{B\subset(C,x)}
\mu^{\mathbb Q}(\F,B,x)
=
\chi(\widehat C)
-
T_{\F}\cdot C
-
\sum_{x\in C\cap\operatorname{Sing}(X)}
\ \sum_{B\subset(C,x)}
\left(1-\frac{1}{|G_x|}\right),
\]
where \(\widehat C\) is the normalization of \(C\), and each
inner sum runs over the local branches of \(C\) at \(x\).
\end{prop}

\begin{proof}
Let
\[
\nu:\widehat C\longrightarrow C
\]
be the normalization. Since \(C\) is \(\F\)-invariant, the
inclusion \(T_{\F}\hookrightarrow T_X\) induces a nonzero
morphism of \(\mathbb Q\)-line bundles
\[
\nu^*T_{\F}\longrightarrow T_{\widehat C}.
\]
We compute its local \(\mathbb Q\)-order at each point of
\(\widehat C\).

Fix \(x\in C\) and a local branch \(B\subset(C,x)\). Suppose
first that \(x\in\operatorname{Sing}(X)\), and let
\[
\pi:(\mathbb C^2,0)\longrightarrow(X,x)
\]
be a local uniformization with group \(G_x\). Write
\[
\pi^{-1}(B)=B_1\cup\cdots\cup B_r,
\]
and let \(H_B\subset G_x\) be the stabilizer of \(B_1\).
The action of \(G_x\) is transitive on these branches, so
\[
r=\frac{|G_x|}{|H_B|}.
\]

Choose a parameter \(t\) on the normalization of \(B_1\) and
a parameter \(s\) on the normalization of \(B\). Since \(H_B\)
acts faithfully on the normalization of \(B_1\), we may choose
these parameters so that
\[
s=t^{|H_B|}.
\]
If \(v\) is a local generator of the lifted foliation
\(\widetilde{\F}\), write
\[
v|_{B_1}=a(t)\frac{\partial}{\partial t},
\qquad
m_B=\operatorname{ord}_0 a(t).
\]
By symmetry, the relative Milnor number on each branch \(B_j\)
is \(m_B\). For the reduced curve \(\pi^{-1}(B)\), however,
the relative Milnor number is
\[
\mu(\widetilde{\F},\pi^{-1}(B),0)
=
r m_B-r+1.
\]
Consequently,
\[
\mu^{\mathbb Q}(\F,B,x)
=
\frac{r m_B-r+1}{|G_x|}.
\]

On the other hand,
\[
d\pi\!\left(a(t)\frac{\partial}{\partial t}\right)
=
|H_B|\,t^{|H_B|-1}a(t)
\frac{\partial}{\partial s}.
\]
Thus the local \(\mathbb Q\)-order of the induced morphism at
the point of \(\widehat C\) corresponding to \(B\) is
\[
\frac{m_B+|H_B|-1}{|H_B|}.
\]
Using \(r=|G_x|/|H_B|\), we obtain
\[
\begin{aligned}
\frac{m_B+|H_B|-1}{|H_B|}
&=
\frac{r m_B}{|G_x|}
+1-\frac{r}{|G_x|}\\
&=
\frac{r m_B-r+1}{|G_x|}
+1-\frac{1}{|G_x|}\\
&=
\mu^{\mathbb Q}(\F,B,x)
+1-\frac{1}{|G_x|}.
\end{aligned}
\]

At a smooth point \(x\in X\), the same computation, with
\(G_x=\{1\}\), gives local order
\[
\mu(\F,B,x)=\mu^{\mathbb Q}(\F,B,x).
\]
The relative Milnor contribution vanishes where the foliation
is nonsingular. Therefore, summing the local orders and
comparing degrees yields
\[
\begin{aligned}
\deg T_{\widehat C}-\deg\nu^*T_{\F}
={}&
\sum_{x\in C\cap\operatorname{Sing}(\F)}
\ \sum_{B\subset(C,x)}
\mu^{\mathbb Q}(\F,B,x)\\
&+
\sum_{x\in C\cap\operatorname{Sing}(X)}
\ \sum_{B\subset(C,x)}
\left(1-\frac{1}{|G_x|}\right).
\end{aligned}
\]
Finally,
\[
\deg T_{\widehat C}=\chi(\widehat C)
\qquad\text{and}\qquad
\deg\nu^*T_{\F}=T_{\F}\cdot C.
\]
Rearranging proves the formula.
\end{proof}

The correction term
\[
\sum_{x\in C\cap\sing(X)}
\left(1-\frac{1}{|G_x|}\right)
\]
measures the discrepancy between the orbifold Euler characteristic of $C$ and
the ordinary Euler characteristic of its normalization; it vanishes exactly
when $C$ avoids the singular locus of $X$. Proposition~\ref{P:Q-milnorrelativo}
is thus the orbifold analogue of the classical Brunella--Khanedani--Suwa
formula, to which it reduces when $X$ is smooth.

\begin{example}
Consider the weighted projective plane
\[
X=\mathbb P(1,1,2)
\]
with quasihomogeneous coordinates $[x:y:z]$ of weights
$1,1,2$, respectively. Let $\F$ be the foliation induced by
the vector field
\[
v=x\frac{\partial}{\partial x}
  +2y\frac{\partial}{\partial y}.
\]
More precisely, this vector field defines a section of $T_X$;
in quasihomogeneous coordinates, it is considered modulo the
weighted Euler vector field
\[
R=x\frac{\partial}{\partial x}
  +y\frac{\partial}{\partial y}
  +2z\frac{\partial}{\partial z}.
\]

Consider the invariant curve
\[
C=\{x=0\}\subset X.
\]
The surface $X$ has a unique singular point,
\[
p=[0:0:1],
\]
whose local uniformizing group is $\mu_2$. The curve $C$
passes through $p$ and is isomorphic, as a curve, to
$\mathbb P^1$. Thus,
\[
\chi(\widehat C)=2.
\]

The vector field $v$ has isolated zeros on $X$. Hence, it
generates the tangent sheaf of the foliation outside a subset
of codimension two and, by saturation,
\[
T_{\F}\simeq\mathcal O_X.
\]
In particular,
\[
T_{\F}\cdot C=0.
\]

We now compute the relative Milnor numbers. At the singular
point $p$, consider the uniformizing chart
\[
\pi:(\mathbb C^2,0)\longrightarrow(X,p),
\qquad
(u,w)\longmapsto [u:w:1],
\]
where $\mu_2$ acts by
\[
(u,w)\longmapsto(-u,-w).
\]
The lifted foliation is defined by
\[
\widetilde v
=
u\frac{\partial}{\partial u}
+
2w\frac{\partial}{\partial w},
\]
and the preimage of the unique branch $B$ of $C$ at $p$ is
\[
\pi^{-1}(B)=\{u=0\}.
\]
Restricting the vector field to this line, we obtain
\[
\widetilde v|_{\{u=0\}}
=
2w\frac{\partial}{\partial w}.
\]
Therefore,
\[
\mu(\widetilde{\F},\pi^{-1}(B),0)=1
\qquad\text{and}\qquad
\mu^{\mathbb Q}(\F,B,p)=\frac12.
\]

Next, consider the smooth point
\[
q=[0:1:0]\in C.
\]
On the chart $y\neq0$, we use the coordinates
\[
a=\frac{x}{y},
\qquad
b=\frac{z}{y^2}.
\]
In these coordinates, the vector field defining $\F$ is
\[
v=-a\frac{\partial}{\partial a}
  -4b\frac{\partial}{\partial b},
\]
and $C=\{a=0\}$. Hence,
\[
v|_C=-4b\frac{\partial}{\partial b}
\qquad\text{and}\qquad
\mu^{\mathbb Q}(\F,C,q)=\mu(\F,C,q)=1.
\]

The only singular points of $\F$ on $C$ are $p$ and $q$.
Consequently, the left-hand side of the formula is
\[
\sum_{x\in C\cap\operatorname{Sing}(\F)}
\ \sum_{B\subset(C,x)}
\mu^{\mathbb Q}(\F,B,x)
=
\frac12+1
=
\frac32.
\]

On the other hand, $C$ has a unique branch at $p$ and
$|G_p|=2$. Therefore, the right-hand side is
\[
\begin{aligned}
\chi(\widehat C)-T_{\F}\cdot C
-
\sum_{x\in C\cap\operatorname{Sing}(X)}
\ \sum_{B\subset(C,x)}
\left(1-\frac{1}{|G_x|}\right)
&=
2-0-\left(1-\frac12\right)\\
&=\frac32.
\end{aligned}
\]
Thus, both sides of the formula agree.
\end{example}

We close the section with the degree bound that drives every argument in what
follows. Recall that a singular point $x$ of $\F$ is \emph{non-dicritical} if
$\F$ admits only finitely many separatrices through $x$.

\begin{prop}\label{P:bounddegree}
Let $X$ be an orbifold projective surface with
$\Cl(X)\cong\mathbb Z$, let $H$ be the ample generator of
$\Cl(X)$, let $\F$ be a foliation on $X$, and let
$C\subset X$ be an $\F$-invariant curve.
Write
$
N_{\F}\equiv dH$ and $
C\equiv mH$,
where $d=\deg(N_{\F})$ and $m=\deg(C)$.
If every singularity of $\F$ lying on $C$ is non-dicritical, then
$
\deg(C)\le \deg(N_{\F}).$

\end{prop}

\begin{proof}
Let $x\in \sing(\F)\cap C$ and let $p:V\longrightarrow U$ be a local
uniformization at $x$. Since $x$ is a non-dicritical singularity of $\F$, every
point $y\in p^{-1}(x)$ is a non-dicritical singularity of $\widetilde{\F}$, so
that $\GSV(\widetilde{\F},\widetilde C,y)\ge0$ by
\cite[Proposition~6]{Brunella1997}. Hence $\QGSV(\F,C,x)\ge0$, and summing over
all singular points of $\F$ on $C$ we obtain $\QGSV(\F,C)\ge0$.
Applying the orbifold GSV index formula, Proposition~\ref{P:Q-GVS-CS},
we obtain
\[
m(d-m)H^2
=
dmH^2-m^2H^2
=
\QGSV(\F,C)
\ge0.
\]
Since $H$ is ample, the Nakai--Moishezon criterion yields $0<H\cdot C=mH^2$, so $m>0$. As $H^2>0$, the above inequality implies $d-m\ge0$, that is, $d\ge m$. This completes the proof.
\end{proof}

\begin{lemma}\label{L:c2smooth}
Let $X$ be a smooth projective surface such that $\operatorname{Pic}(X)\simeq\mathbb{Z}$. Then $c_2(X)\geq 3$.
\end{lemma}

\begin{proof}
Since $\operatorname{Pic}(X)\simeq \mathbb{Z}$, the group
$\operatorname{Pic}^0(X)$ is trivial. Hence $q(X)=h^1(X,\mathcal{O}_X)=0$, and therefore $b_1(X)=2q(X)=0$. Moreover, $\rho(X)=\operatorname{rank}\operatorname{NS}(X)=1$, and since $\rho(X)\leq b_2(X)$, we obtain $b_2(X)\geq 1$. Finally, by the topological Gauss--Bonnet formula, $c_2(X)=e(X)=2-2b_1(X)+b_2(X)$; thus $c_2(X)=2+b_2(X)\geq 3$.
\end{proof}

Let $\mathcal F$ be a foliation on an orbifold surface $X$, and let
$x\in X$ be a singular point of $\mathcal F$. We define the algebraic
multiplicity $m_x(\mathcal F)$ as follows.

If $x$ is a regular point of $X$, let $v$ be a holomorphic vector field
defining $\mathcal F$ locally around $x$. Writing
\[
v=v_r+v_{r+1}+\cdots,
\]
where $v_j$ is homogeneous of degree $j$ and $v_r\neq0$, we set
\[
m_x(\mathcal F)=r.
\]

If $x$ is a singular point of $X$, we lift $\mathcal F$ to a local
uniformizing chart $\widetilde U\subset\mathbb C^2$ around a lift of
$x$. If $v$ is a holomorphic vector field defining the lifted foliation,
we define $m_x(\mathcal F)$ as above, namely as the degree of the first
nonzero homogeneous term in the expansion of $v$.
Similarly, we say that $x$ is a \emph{nilpotent singularity}
(respectively, a \emph{saddle-node}) of $\mathcal F$ if
the lifted foliation $\widetilde{\mathcal F}$ has a
nilpotent singularity (respectively, a saddle-node)
at some point of $\pi^{-1}(x)$.
In the following section, we study the arithmetic obstructions imposed by the condition $m_p(\mathcal F)=1$ for foliations on weighted projective planes.

\section{Foliations with a unique singularity}

We now consider foliations with a unique singularity on a weighted
projective plane $\mathbb P(a,b,c)$. In particular, we assume that the
foliation has a unique singularity $p$ and that its algebraic
multiplicity is equal to one.

Under these assumptions, the local description of foliations on weighted
projective planes provides arithmetic conditions on the normal degree
$d$. More precisely, the congruence of $d$ modulo each weight
$a$, $b$, and $c$ is determined by the position of the unique
singularity $p$ with respect to the coordinate points of
$\mathbb P(a,b,c)$. When $p$ is different from a coordinate point, the
foliation is regular at that point, whereas if $p$ coincides with a
coordinate point, the assumption that the algebraic multiplicity is one
gives three possible congruence classes.

The following lemma collects these arithmetic conditions.

\begin{lemma} \label{L:localobstruction}
Let $\mathcal F$ be a foliation of normal degree $d$ on the weighted
projective plane
\[
\mathbb P(a,b,c),
\]
where $a,b,c$ are pairwise coprime. Assume that $p$ is the unique
singularity of $\mathcal F$ and that $p$ has ordinary multiplicity one.
Denote the coordinate points by
\[
p_0=(1:0:0),\qquad
p_1=(0:1:0),\qquad
p_2=(0:0:1).
\]

Then the following possibilities hold.

\begin{enumerate}
    \item If $p\notin\{p_0,p_1,p_2\}$, then
    \[
    \boxed{
    \begin{aligned}
    d&\equiv b\ \text{or}\ c &&\pmod a,\\
    d&\equiv a\ \text{or}\ c &&\pmod b,\\
    d&\equiv a\ \text{or}\ b &&\pmod c.
    \end{aligned}}
    \]

    \item If $p=p_0$, then
    \[
    \boxed{d\equiv 2b,\ b+c,\ 2c\pmod a,}
    \]
    and
    \[
    d\equiv a\ \text{or}\ c\pmod b,\qquad
    d\equiv a\ \text{or}\ b\pmod c.
    \]

    \item If $p=p_1$, then
    \[
    \boxed{d\equiv 2a,\ a+c,\ 2c\pmod b,}
    \]
    and
    \[
    d\equiv b\ \text{or}\ c\pmod a,\qquad
    d\equiv a\ \text{or}\ b\pmod c.
    \]

    \item If $p=p_2$, then
    \[
    \boxed{d\equiv 2a,\ a+b,\ 2b\pmod c,}
    \]
    and
    \[
    d\equiv b\ \text{or}\ c\pmod a,\qquad
    d\equiv a\ \text{or}\ c\pmod b.
    \]
\end{enumerate}
\end{lemma}

\begin{proof}
We use the local description of foliations on weighted projective planes
established in \cite[Equation (2.7), page 5]{lizarbe2023density}.
In a local chart
\[
U_0\simeq\mathbb C^2/\mu_a,
\]
the defining $1$-form decomposes as
\[
\omega=\omega_r+\omega_{r+a}+\cdots+\omega_d,
\]
where the nonzero components are quasi-homogeneous and their degrees
differ by $a$. Hence
\[
\boxed{d\equiv r\pmod a.}
\]

At $p_1$ and $p_2$, the analogous decompositions give
\[
d\equiv r\pmod b
\qquad\text{and}\qquad
d\equiv r\pmod c,
\]
respectively. Moreover, $r$ is precisely the algebraic multiplicity at
the corresponding coordinate point.

If $p\neq p_0,p_1,p_2$, then all three coordinate points are regular for
$\mathcal F$, since $p$ is the unique singularity of the foliation. At
$p_0$, the two local coordinates have weights $b$ and $c$; hence the
first nonzero term has degree $b$ or $c$. Therefore
\[
d\equiv b\ \text{or}\ c\pmod a.
\]
Similarly,
\[
d\equiv a\ \text{or}\ c\pmod b,
\qquad
d\equiv a\ \text{or}\ b\pmod c.
\]
This proves (1).

If $p=p_0$, multiplicity one means that the first nonzero terms are
obtained from the linear part in the two local coordinates, whose
weights are $b$ and $c$. Their possible weighted degrees are
\[
2b,\qquad b+c,\qquad 2c.
\]
Consequently,
\[
d\equiv2b,\ b+c,\ 2c\pmod a.
\]
Since $p_1$ and $p_2$ are regular for $\mathcal F$, we also have
\[
d\equiv a\ \text{or}\ c\pmod b,
\qquad
d\equiv a\ \text{or}\ b\pmod c.
\]
This proves (2). The cases $p=p_1$ and $p=p_2$ follow in exactly the
same way.
\end{proof}

We now apply the previous lemma to some special weighted projective
planes, where the resulting systems of congruences can be explicitly
solved using the Chinese Remainder Theorem.

\begin{lemma} \label{L:particularlocalobstruction}
The congruence conditions of the previous lemma can be explicitly
determined in the following cases.

\begin{enumerate}
    \item For $\mathbb P(1,1,k)$, let
    \[
    p_2=(0:0:1).
    \]
    Then
    \[
    p\neq p_2
    \quad\Longrightarrow\quad
    d\equiv1\pmod k,
    \]
    while
    \[
    p=p_2
    \quad\Longrightarrow\quad
    d\equiv2\pmod k.
    \]

    \item For $\mathbb P(1,a,b)$, with $\gcd(a,b)=1$, if
    \[
    p\notin\{p_1,p_2\},
    \]
    then
    \[
    d\equiv1\ \text{or}\ b\pmod a,
    \qquad
    d\equiv1\ \text{or}\ a\pmod b.
    \]

    If $p=p_1$, then
    \[
    d\equiv2,\ b+1,\ 2b\pmod a,
    \qquad
    d\equiv1\ \text{or}\ a\pmod b.
    \]

    If $p=p_2$, then
    \[
    d\equiv1\ \text{or}\ b\pmod a,
    \qquad
    d\equiv2,\ a+1,\ 2a\pmod b.
    \]

    \item For $\mathbb P(2,3,5)$, the complete list of congruence
    classes modulo $30$ is
    \[
    \begin{array}{c|c}
    \text{position of }p & d\pmod{30}\\
    \hline
    p\notin\{p_0,p_1,p_2\} & 17,\ 23\\
    p=p_0 & 2,\ 8\\
    p=p_1 & 7,\ 13\\
    p=p_2 & 5,\ 11,\ 29.
    \end{array}
    \]
\end{enumerate}
\end{lemma}

\begin{proof}
The assertions follow by substituting the corresponding weights into
the congruences of the previous lemma. In the cases
$\mathbb P(1,a,b)$ and $\mathbb P(2,3,5)$, the resulting systems involve
pairwise coprime moduli and can therefore be solved using the Chinese
Remainder Theorem.

For $\mathbb P(2,3,5)$, for instance, if
\[
p\notin\{p_0,p_1,p_2\},
\]
the previous lemma gives
\[
d\equiv1\pmod2,\qquad
d\equiv2\pmod3,\qquad
d\equiv2\ \text{or}\ 3\pmod5.
\]
The Chinese Remainder Theorem yields
\[
d\equiv17\quad\text{or}\quad23\pmod{30}.
\]
The other cases are obtained in the same way.
\end{proof}

\section{Foliations on Orbifold Projective Surfaces with Cyclic Divisor Class Group}\label{S:classone}

Let $X$ be an orbifold projective surface whose divisor class group is
cyclic and which satisfies a suitable topological condition. In this section
we prove that a foliation on $X$ with a unique singularity and sufficiently
large normal degree admits no invariant algebraic curve.

The argument splits according to the type of the singularity. If its
algebraic multiplicity is one, the singularity is either a saddle-node or
nilpotent, and we treat the two cases separately, in
Propositions~\ref{P:orbifoldsaddle-node} and~\ref{P:orbifoldnilpotent};
when $X$ is smooth or toric the nilpotent bound can be sharpened
(Proposition~\ref{P:smooth-toricnilpotent}). Combining the two cases yields
Theorem~\ref{T:classone}. We begin with a lemma that bounds the orbifold GSV
index of an invariant curve.

\begin{lemma}\label{L:classone}
Let $X$ be an orbifold projective surface with $\Cl(X)=\mathbb{Z}$, let $\F$
be a foliation on $X$, and let $C\subset X$ be a complete $\F$-invariant
curve. Then:
\begin{enumerate}
\item
$
\QGSV(\F,C)\le \frac{N_{\F}^{2}}{4},
$ 
with equality if and only if $\deg(C)=\deg(N_{\F})/2$.

\item If every singularity of $\F$ lying on $C$ is non-dicritical and
$\deg(K_{\F})\ge0$, then
\[
C\cdot K_{\F}\le K_{\F}\cdot N_{\F}.
\]
\end{enumerate}
\end{lemma}

\begin{proof}
Let $H$ be the ample generator of $\Cl(X)$. Since $\Cl(X)=\mathbb Z$, there exist positive integers $m$ and $d$ such that $C\equiv mH$ and $N_{\mathcal F}\equiv dH$.
Since $H$ is ample, we have $H^2\in\mathbb Q_{>0}$.

For (1), the GSV-index formula yields
\[
\QGSV(\mathcal F,C)
=
N_{\mathcal F}\cdot C-C^2.
\]
Hence $\QGSV(\mathcal F,C)=(md-m^2)H^2$.
Since
\[
md-m^2
=
-\left(m-\frac d2\right)^2+\frac{d^2}{4}
\le
\frac{d^2}{4},
\]
and $H^2>0$, we obtain
\[
\GSV(\mathcal F,C)
\le
\frac{d^2H^2}{4}
=
\frac{N_{\mathcal F}^2}{4}.
\]

For (2), by Proposition~\ref{P:bounddegree}, $\deg(C)\le\deg(N_{\mathcal F})$. Since $\deg(K_{\mathcal F})\ge0$, we have
\[
\deg(C)\deg(K_{\mathcal F})
\le
\deg(N_{\mathcal F})\deg(K_{\mathcal F}).
\]
Moreover, $K_{\mathcal F}\equiv kH$ for some rational number $k\ge0$, because $\Cl(X)=\mathbb Z$.
Therefore,
\[
C\cdot K_{\mathcal F}
=
mkH^2,
\qquad
K_{\mathcal F}\cdot N_{\mathcal F}
=
dkH^2.
\]
Since
\[
\deg(C)=mH^2,\qquad
\deg(N_{\mathcal F})=dH^2,
\]
the inequality $\deg(C)\le\deg(N_{\mathcal F})$ is equivalent to $m\le d$. As $kH^2\ge0$, it follows that $mkH^2\le dkH^2$, that is,
\[
C\cdot K_{\mathcal F}
\le
K_{\mathcal F}\cdot N_{\mathcal F}.
\]
This completes the proof.
\end{proof}

We first treat the saddle-node case.

\begin{prop}\label{P:orbifoldsaddle-node}
Let $X$ be an orbifold projective surface with
$\Cl(X)\cong\mathbb Z$, let $H$ be the ample generator of
$\Cl(X)$, and let $\F$ be a foliation on $X$ with a unique singularity
$x$. Write $N_{\F}\equiv dH$, where $d>0$. If $x$ is a saddle-node
singularity and $d>\max\{0,B_1\}$, where
\[
B_1=
\frac{2}{3H^2}
\left(
c_1(T_X)\cdot H+
\sqrt{\bigl(c_1(T_X)\cdot H\bigr)^2
-3H^2\bigl(c_2(X)-1\bigr)}
\right),
\]
then $\F$ admits no invariant algebraic curve.
\end{prop}

\begin{proof}
Since $N_{\F}\equiv dH$, we have
\[
\mu^{\Q}(\F)
=
H^2d^2-c_1(T_X)\cdot Hd+c_2(X).
\]

Assume, by contradiction, that $\F$ admits an invariant algebraic curve
$C\subset X$.
By the orbifold GSV-index formula, $\QGSV(\F,C)=N_{\F}\cdot C-C^2$.
Hence, by Lemma~\ref{L:classone},
\[
\QGSV(\F,C)
\le
\frac{N_{\F}^2}{4}
=
\frac{H^2d^2}{4}.
\]

Let $p:V\longrightarrow U$ be a local orbifold chart around $x$, let $\widetilde{\F}=p^{*}\F$, and let $\widetilde C=p^{-1}(C)$.
Since $x$ is a saddle-node singularity of $\F$, every point of
$p^{-1}(x)$ is a saddle-node singularity of $\widetilde{\F}$.
Therefore, by~\cite[Proposition~2]{FernandezPuchuriRosas2021},
\[
\GSV(\widetilde{\F},\widetilde C,y)
\ge
\mu(\widetilde{\F},y)-1
\]
for every $y\in p^{-1}(x)$.
Averaging over the fiber, we obtain $\QGSV(\F,C,x)\ge\mu^{\Q}(\F,x)-1$.
Since $x$ is the unique singularity of $\F$,
\[
\QGSV(\F,C)=\QGSV(\F,C,x)
\quad\text{and}\quad
\mu^{\Q}(\F)=\mu^{\Q}(\F,x).
\]
Hence
\[
\frac{H^2d^2}{4}
\ge
H^2d^2-c_1(T_X)\cdot Hd+c_2(X)-1,
\]
or equivalently,
\[
\frac34H^2d^2
-\bigl(c_1(T_X)\cdot H\bigr)d
+c_2(X)-1
\le0.
\]
If the discriminant is negative, then the left-hand side is positive
for every $d$, a contradiction. Thus the discriminant is non-negative,
and solving the quadratic inequality gives
\[
d\le
\frac{2}{3H^2}
\left(
c_1(T_X)\cdot H+
\sqrt{
\bigl(c_1(T_X)\cdot H\bigr)^2
-3H^2\bigl(c_2(X)-1\bigr)}
\right)
=B_1,
\]
contradicting the assumption that $d>\max\{0,B_1\}$.
Therefore, $\F$ admits no invariant algebraic curve.
\end{proof}

\begin{example}\label{E:weightedplane(1,1,k)-saddle-node}
Consider the weighted projective plane $X=\mathbb P(1,1,k)$, $k>1$.
Let $H$ denote the ample generator of $\Cl(X)$. Then
\[
H^2=\frac1k,\qquad
c_1(T_X)\cdot H=\frac{k+2}{k},\qquad
c_2(X)=2+\frac1k.
\]
Hence
\[
B_1
=
\frac23\left(k+2+\sqrt{k^2+k+1}\right).
\]
Since $k^2+k+1<(k+1)^2$, we obtain
\[
B_1
<
\frac23(2k+3)
=
\frac{4k}{3}+2
<
2k+1.
\]
Therefore, the hypothesis of Proposition~\ref{P:orbifoldsaddle-node} is
satisfied whenever $d\ge 2k+1$.
Consequently, every foliation on $\mathbb P(1,1,k)$ with a unique
saddle-node singularity and normal class $N_{\mathcal F}\equiv dH$, $d\ge 2k+1$, admits no invariant algebraic curve.
\end{example}

\begin{example}\label{E:weightedplane(1,2,3)-saddle-node}
Consider the weighted projective plane $X=\mathbb P(1,2,3)$.
Let $H$ be the ample generator of $\Cl(X)$. Then
\[
H^2=\frac16,\qquad
c_1(T_X)\cdot H=1,\qquad
c_2(X)=\frac{11}{6}.
\]
Hence
\[
B_1
=
\frac{2}{3H^2}
\left(
c_1(T_X)\cdot H+
\sqrt{\bigl(c_1(T_X)\cdot H\bigr)^2
-3H^2\bigl(c_2(X)-1\bigr)}
\right)
=
4+2\sqrt2.
\]
Since $4+2\sqrt2<7$, the hypothesis of Proposition~\ref{P:orbifoldsaddle-node} is satisfied whenever $d\ge7$.
Consequently, every foliation on $\mathbb P(1,2,3)$ with a unique
saddle-node singularity and normal class $N_{\mathcal F}\equiv dH$, $d\ge7$, admits no invariant algebraic curve.
\end{example}

We turn to the nilpotent case.

\begin{prop}\label{P:orbifoldnilpotent}
Let $X$ be an orbifold projective surface with
$\Cl(X)\cong\mathbb Z$ and $c_2(X)>1$. Let $H$ be the ample generator of
$\Cl(X)$, and let $\F$ be a foliation on $X$ with a unique nilpotent
singularity $x$. Write $N_{\F}\equiv dH$, where $d>0$. Assume that
\[
d>\max\left\{
0,\,
\frac{c_1(T_X)\cdot H}{H^2},\,
B_2
\right\},
\]
where
\[
B_2=
\frac{
c_1(T_X)\cdot H+
\sqrt{\bigl(c_1(T_X)\cdot H\bigr)^2
-2H^2\left(c_2(X)+\frac{1}{|G_x|}\right)}
}{H^2},
\]
and $G_x$ denotes the local orbifold group at $x$.
Then $\F$ admits no invariant algebraic curve.
\end{prop}

\begin{proof}
Since $N_{\F}\equiv dH$, we have
\[
\mu^{\Q}(\F)
=
H^2d^2-c_1(T_X)\cdot Hd+c_2(X).
\]

Assume, by contradiction, that there exists an invariant algebraic curve
$C\subset X$.
By the orbifold GSV-index formula, $\QGSV(\F,C)=N_{\F}\cdot C-C^2$.
Hence, by Lemma~\ref{L:classone},
\[
\QGSV(\F,C)
\le
\frac{N_{\F}^2}{4}.
\]    

    Let $p:V\longrightarrow U$ be a local orbifold chart around $x$, let $\widetilde{\F}=p^{*}\F$, and let $\widetilde C=p^{-1}(C)$.
Since $x$ is a nilpotent singularity of $\F$, every point of
$p^{-1}(x)$ is a nilpotent singularity of $\widetilde{\F}$.
Therefore,~\cite[Proposition~3]{FernandezPuchuriRosas2021} implies that one of the following alternatives holds for every $y\in p^{-1}(x)$:

\begin{enumerate}
\item $BB(\widetilde{\F},y)=0$;

\item
$
BB(\widetilde\F,y)\ge2\bigl(\mu(\widetilde\F,y)+1\bigr);
$

\item$
\GSV(\widetilde\F,\widetilde C,y)
\ge
\frac{\mu(\widetilde\F,y)+1}{2};$

\item
\(\widetilde{C}\) is singular only at points \(y\) where it has exactly two local branches; moreover, every such point \(y\) is a non-dicritical singularity of \(\widetilde{\mathcal F}\) and satisfies
\[
\sum_{\widetilde{B}\subset(\widetilde C,y)}
\mu(\widetilde{\mathcal F},\widetilde B,y)
=
\mu(\widetilde{\mathcal F},y)+1.
\]
\end{enumerate}

Averaging over the fiber, we obtain the corresponding alternatives for $\F$:

\begin{enumerate}
\item $BB^{\Q}(\F,x)=0$;

\item
$
BB^{\Q}(\F,x)
\ge
2\mu^{\Q}(\F,x)+\frac{2}{|G_x|};
$

\item
$
\QGSV(\F,C,x)
\ge
\frac{\mu^{\Q}(\F,x)}{2}
+\frac{1}{2|G_x|};
$

\item Note that the two branches are not conjugate under the action of the local orbifold group at $x$. Indeed, if they were, then the orbifold Camacho--Sad indices of the two local branches at $y$ would coincide. By \cite[Case 2a, p.~7]{FernandezPuchuriRosas2021}, this would imply that $BB(\widetilde{\mathcal F},y)=0$, and consequently $BB(\mathcal F,x)=0$, a contradiction. 
Therefore, $C$ is singular only at the point $x$, where it has exactly two local branches; moreover, $x$ is a non-dicritical singularity of $\F$ and satisfies
\[
\sum_{B\subset(C,x)}
\mu^{\Q}(\F,B,x)
=
\mu^{\Q}(\F,x)+\frac{1}{|G_x|}.
\]

\end{enumerate}

Since $x$ is the unique singularity of $\F$, we have
\[
BB^{\Q}(\F)=BB^{\Q}(\F,x),
\quad
\mu^{\Q}(\F)=\mu^{\Q}(\F,x),\quad
\QGSV(\F,C)=\QGSV(\F,C,x).\]

If {\rm(1)} holds, then the Baum--Bott formula yields $0=BB^{\Q}(\F)=N_{\F}^2=d^2H^2$, contradicting the assumption that $d>0$.

Suppose that {\rm(2)} holds. Since $BB^{\Q}(\F)=N_{\F}^2$, we obtain
\[
N_{\F}^2\ge 2\mu^{\Q}(\F)+\frac{2}{|G_x|}.
\]

Suppose now that {\rm(3)} holds. Then
\[
\frac{N_{\F}^2}{4}
\ge
\QGSV(\F,C)\ge
\frac{\mu^{\Q}(\F)}{2}+\frac{1}{2|G_x|},
\]
which again implies
\[
N_{\F}^2\ge 2\mu^{\Q}(\F)+\frac{2}{|G_x|}.
\]
Using the expression
\[
\mu^{\Q}(\F)
=
H^2d^2-c_1(T_X)\cdot Hd+c_2(X),
\]
we obtain
\[
H^2d^2
-
2(c_1(T_X)\cdot H)d
+
2\left(c_2(X)+\frac{1}{|G_x|}\right)
\le 0.
\]
If the discriminant of this quadratic polynomial is negative, then,
since its leading coefficient is positive, the polynomial is positive
for every $d$, a contradiction. Hence the discriminant is
non-negative, and solving the quadratic inequality gives
\[
d
\le
\frac{
c_1(T_X)\cdot H+
\sqrt{(c_1(T_X)\cdot H)^2
-2H^2\left(c_2(X)+\frac{1}{|G_x|}\right)}
}{H^2}
=
B_2,
\]
contradicting
\[
d>
\max\left\{
0,
\frac{c_1(T_X)\cdot H}{H^2},
B_2
\right\}.
\]

Finally, suppose that {\rm(4)} holds. By
Proposition~\ref{P:Q-milnorrelativo},
\[
\chi(\widehat C)
=
\mu^{\Q}(\F)
+\frac{1}{|G_x|}
-C\cdot K_{\F}
+\sum_{y\in C\cap\sing(X)}
\left(1-\frac1{|G_y|}\right),
\]
where $\widehat C$ denotes the normalization of $C$. By
Lemma~\ref{L:classone}, $C\cdot K_{\F}\le K_{\F}\cdot N_{\F}$, and therefore
\[
\chi(\widehat C)
\ge
\mu^{\Q}(\F)
+1
-N_{\F}\cdot K_{\F}
+\sum_{y\in C\cap\sing(X)\setminus\{x\}}
\left(1-\frac1{|G_y|}\right).
\]
Since $K_{\F}=N_{\F}-c_1(T_X)=-T_{\F}$, we obtain
\[
\chi(\widehat C)
\ge
\mu^{\Q}(\F)
+1
+N_{\F}\cdot T_{\F}
=
c_2(X)+1.
\]
Since $c_2(X)>1$, we have $\chi(\widehat C)>2$,
which is impossible because the normalization of an irreducible
projective curve is a smooth projective curve whose Euler
characteristic is at most $2$.
Thus every possible case leads to a contradiction. Therefore,
$\F$ admits no invariant algebraic curve.
\end{proof}

When $X$ is smooth or toric, the bound of
Proposition~\ref{P:orbifoldnilpotent} can be sharpened: the strict
inequality $d>B_2$ may be replaced by $d\ge B_2$.

\begin{prop}\label{P:smooth-toricnilpotent}
Under the hypotheses of Proposition~\ref{P:orbifoldnilpotent}, suppose furthermore that $X$ is smooth or toric. If
\[
d>\max\left\{0,\frac{c_1(T_X)\cdot H}{H^2}\right\}
\quad\text{and}\quad
d\ge B_2,
\]
then $\F$ admits no invariant algebraic curve.
\end{prop}

\begin{proof}
Suppose, by contradiction, that there exists an algebraic curve
$C\subset X$ invariant under $\F$. By the proof of
Proposition~\ref{P:orbifoldnilpotent}, it suffices to consider the case
where $C$ is quasi-smooth at $x$,
\[
BB^{\Q}(\F,x)=2\mu^{\Q}(\F,x)+\frac{2}{|G_x|}
\]
and
\[
\QGSV(\F,C,x)
=\frac{\mu^{\Q}(\F,x)}{2}
+\frac{1}{2|G_x|}.
\]
Since $x$ is the unique singularity of $\F$, we have
$
\QGSV(\F,C)=\frac{N_{\F}^{2}}{4}.
$
Moreover, Lemma~\ref{L:classone} implies that
$
\deg(C)=\frac{\deg(N_{\F})}{2}.
$
On the other hand, Proposition~\ref{P:Q-milnorrelativo} yields
\[
\chi(\widetilde C)
=
\frac{\mu(\F)}{2}
+\frac{1}{2|G_x|}
-K_{\F}\cdot\frac{N_{\F}}{2}
+\left(1-\frac1{|G_x|}\right),
\]
and therefore
\[
\chi(\widetilde C)
=
\frac{c_2(X)}{2}
+1-\frac{1}{2|G_x|}.
\]
Since $\chi(\widetilde C)\in\mathbb Z$ and $c_2(X)>1$, this is impossible unless
\begin{equation}\label{E:c_2}
c_2(X)=2+\frac1{|G_x|}
\quad\text{and}\quad
\chi(\widetilde C)=2.
\end{equation}
Assume first that $X$ is smooth. Then $C$ is smooth and 
\[
2=\chi(C)
=\GSV(\F,C)-K_{\F}\cdot C
=\frac{N_{\F}^{2}}{4}
-K_{\F}\cdot\frac{N_{\F}}{2}.
\]
Writing $K_{\F}\equiv kH$, we obtain $8=d(d-2k)H^2$. Since
$
d>\frac{c_1(T_X)\cdot H}{H^2},
$
we have
\[
k=d-\frac{c_1(T_X)\cdot H}{H^2}>0.
\]
Hence $d-2k>0$, so both $d$ and $d-2k$ are positive integers. Since $d(d-2k)H^2=8$, it follows that $H^2$ divides $8$. Examining all positive factorizations
of $8/H^2$, we conclude that the only possibility is $H^2=1$, $d=4$ and $k=1$. Consequently, $K_X=K_{\F}\otimes N_{\F}^{*}=-3H$, and therefore $c_1(T_X)=3H$. It follows from the Kobayashi--Ochiai theorem~\cite[Corollary, p.~32]{KobayashiOchiai1973} that $X\simeq\mathbb P^2$. Moreover, $\deg(C)=2$. Hence $C$ is a smooth conic. A direct computation then shows that
$\F$ admits an invariant line, a contradiction.

In the toric case, it follows from \cite[Lemma 2.11]{BatyrevCox1994} and
(\ref{E:c_2}) that
\[
X\simeq \Pj(1,1,k),\qquad
H^2=\frac{1}{k},\qquad
c_1(T_X)\cdot H=\frac{k+2}{k},\qquad
c_2(X)=\frac{2k+1}{k},
\]
and the local orbifold group of $X$ at the singular point has order $|G_x|=k$.
Furthermore,
\[
2=\chi(C)
=\GSV(\F,C)-K_{\F}\cdot C
=\frac{N_{\F}^{2}}{4}
-K_{\F}\cdot\frac{N_{\F}}{2}
+1-\frac{1}{|G_x|}.
\]
Hence,
\[
8=d\bigl(d-2(d-(k+2))\bigr)H^2
+4\left(1-\frac{1}{|G_x|}\right).
\]
Substituting $H^2=\frac1k$ and $|G_x|=k$, we obtain $d^2-(2k+4)d+4(2k+1)=0$. This equation has two solutions, $d=2k+2$ and $d=2$.
Since
\[
d>\frac{c_1(T_X)\cdot H}{H^2}=k+2,
\]
only the first solution is admissible. Therefore $d=2k+2$ and consequently $\deg(C)=k+1$. Hence $C$ is quasi-smooth. A direct computation then shows that
$\F$ admits an invariant line, a contradiction.
\end{proof}

\begin{example}\label{E:nilpotent(1,1,k)}
Let $X=\mathbb P(1,1,k)$ with $k\ge1$, and let $H$ denote the ample
generator of $\Cl(X)$. Let $\F$ be a foliation on $X$ with a unique
nilpotent singularity such that $N_{\F}\equiv dH$.
Since
\[
H^2=\frac1k,\qquad
c_1(T_X)=(k+2)H,\qquad
c_2(X)=2+\frac1k,
\]
we have
$
\frac{c_1(T_X)\cdot H}{H^2}=k+2.
$
Moreover,
$
c_2(X)+\frac1{|G_x|}
=
2+\frac2k,
$
and therefore
\[
B_2
=
\frac{
\frac{k+2}{k}
+
\sqrt{
\left(\frac{k+2}{k}\right)^2
-
\frac{2}{k}\left(2+\frac2k\right)
}
}{\frac1k}
=
2k+2.
\]
Hence,
\[
d\ge B_2
\quad\Longleftrightarrow\quad
d\ge2k+2.
\]
Therefore, by Proposition~\ref{P:smooth-toricnilpotent}, if $d\ge2k+2$, then $\F$ admits no invariant algebraic curve.
\end{example}

Combining the two cases, we obtain the main result of this section.

\begin{thm}\label{T:classone}
Let $X$ be an orbifold projective surface with
$\Cl(X)\cong\mathbb Z$ and $c_2(X)>1$. Let $H$ denote the ample generator
of $\Cl(X)$, and let $\F$ be a foliation on $X$ with a unique singularity
$x$. Write $N_{\F}\equiv dH$, where $d=\deg(N_{\F})$. Assume that the algebraic multiplicity of
$\F$ at $x$ is one.
If
\[
d>\max\left\{
0,\,
\frac{c_1(T_X)\cdot H}{H^2},\,
B_0,\,
B_1,\,
B_2
\right\},
\]
where
\[
B_0=
\frac{
c_1(T_X)\cdot H+
\sqrt{\bigl(c_1(T_X)\cdot H\bigr)^2
-4H^2\bigl(c_2(X)-2\bigr)}
}{2H^2},
\]
and $B_1$ and $B_2$ are the constants defined in
Propositions~\ref{P:orbifoldsaddle-node} and
\ref{P:orbifoldnilpotent}, respectively,
then $\F$ admits no invariant algebraic curve.
Moreover, if $X$ is smooth or toric, then the same conclusion holds under
the weaker assumptions
\[
d>\max\left\{0,\frac{c_1(T_X)\cdot H}{H^2},B_0,\,
B_1\,\right\}
\quad\text{and}\quad
d\ge B_2.
\]
\end{thm}

\begin{proof}
Since $x$ is the unique singularity of $\F$, the Baum--Bott formula gives
\[
\mu^{\Q}(\F,x)
=
H^2d^2-c_1(T_X)\cdot Hd+c_2(X).
\]
Since $d>B_0$, it follows that $\mu^{\Q}(\F,x)\ge2$.
Moreover, as the algebraic multiplicity of $x$ is one, the singularity is
either a saddle-node or a nilpotent singularity. Hence, under the first set
of hypotheses, the conclusion follows from
Propositions~\ref{P:orbifoldsaddle-node} and
\ref{P:orbifoldnilpotent}.

If, furthermore, $X$ is smooth or toric, then the conclusion follows from
Propositions~\ref{P:orbifoldsaddle-node} and
\ref{P:smooth-toricnilpotent}.
\end{proof}

\begin{remark}
    The technical condition $c_2(X)>1$ in Proposition \ref{P:orbifoldnilpotent} is not necessary if $X$ is a smooth surface and $\Pic(X)=\mathbb{Z}$, by Lemma \ref{L:c2smooth}.
\end{remark}

\section{Applications}\label{S:applications}

We now specialize Theorem~\ref{T:classone} to weighted projective planes,
where the constants $B_0$, $B_1$ and $B_2$ can be computed explicitly, and we
show that the resulting bound is sharp. The first result extends
\cite[Theorem~1]{FernandezPuchuriRosas2021} from $\Proj^2$ to
$\Pj(1,1,k)$.

%
%
%
%

\begin{thm}\label{T:(1,1,k)}
Let $\mathcal F$ be a foliation on $\mathbb P(1,1,k)$, $k\geq1$, of
normal degree $d$, with a unique singular point $p$ and
$m_p(\mathcal F)=1$. Let $p_2=(0:0:1)$. Then:

\begin{enumerate}
\item If $d\geq2k+2$, then $\mathcal F$ admits no invariant algebraic
curve.

\item If $d=2k+1$, then $\mathcal F$ admits no invariant algebraic
curve if $p$ is a saddle-node, whereas it may admit one if $p$ is
nilpotent.
\end{enumerate}

Moreover,

$$
d\equiv1\pmod{k}\quad\text{if }p\neq p_2,\qquad
d\equiv2\pmod{k}\quad\text{if }p=p_2.
$$

\end{thm}

\begin{proof}
Assertion~(1) follows from Theorem~\ref{T:classone},
Example~\ref{E:nilpotent(1,1,k)}, and
Example~\ref{E:weightedplane(1,1,k)-saddle-node}.

For assertion~(2), if $p$ is a saddle-node singularity, then the conclusion
follows from Example~\ref{E:weightedplane(1,1,k)-saddle-node}. On the other
hand, if $p$ is a nilpotent singularity, the existence of foliations
admitting invariant algebraic curves is established in the following
remark. 

Finally, the congruence statements in the ``Moreover'' part follow from Lemma~\ref{L:particularlocalobstruction}. This completes the proof.
\end{proof}

\begin{remark}\label{R:nilpontentfirstintegral}
For every $k$ there are foliations of normal degree $d=2k+1$ on
$\Pj(1,1,k)$ with a unique nilpotent singularity that do admit an invariant
algebraic curve. They are defined by the $1$-form
\[
\omega=
-kz^2\,dx
+kzy^k\,dy
+(kxz-y^{k+1})\,dz,
\]
which can be written as
\[
\omega
=
kz\,d\!\left(-xz+\frac{1}{k+1}y^{k+1}\right)
-(k+1)\left(-xz+\frac{1}{k+1}y^{k+1}\right)\,dz.
\]
Hence, $\F$ admits the rational first integral
\[
\frac{-xz+\frac{1}{k+1}y^{k+1}}{z^{k+1}}.
\]
\end{remark}

We now classify the foliations of normal degree $2k+1$ on
$\Pj(1,1,k)$, where $k>1$, having a unique singularity.

\begin{thm}\label{T:case2k+1(1,1k,)}
Let $\F$ be a foliation on $\Pj(1,1,k)$ of normal degree $2k+1$ having a unique singular point $p$.
Then exactly one of the following alternatives occurs:

\begin{enumerate}
\item $p=[0:0:1]$, $m_p(\F)=k$, and $\F$ is a Riccati foliation. In particular, $\F$ admits an invariant algebraic curve.

\item After an automorphism of $\Pj(1,1,k)$, we may assume that
$p=[1:0:0]$. In suitable affine coordinates $(y,z)$, the foliation $\F$
is induced by one of the following $1$-forms:

\begin{enumerate}
\item
\[
\omega=
\Bigl(
y^{2k-1}
+z\sum_{j=2}^{k}a_jy^{j-1}
+kz^2
\Bigr)\,dy
-yz\,dz,
\qquad a_j\in\C .
\]

\item
\[
\omega=
\Bigl(
(z-y^{k+1})
\sum_{j=1}^{k-1}b_jy^{j-1}
+kzy^k
\Bigr)\,dy
+(z-y^{k+1})\,dz,
\qquad b_j\in\C .
\]

\item
\[
\omega=
\Bigl(
(z-y^{k+1}-zy)
\sum_{j=1}^{k-1}b_jy^{j-1}
+kz(y^k+z)
\Bigr)\,dy
+(z-y^{k+1}-zy)\,dz,
\qquad b_j\in\C .
\]
\end{enumerate}

\end{enumerate}
\end{thm}

\begin{proof}
If $p=[0:0:1]$, then $m_p(\F)=k$. In the affine chart $z=1$, the foliation is induced by the $1$-form
\[
\omega_0=
\left(A_{k-1}-yH_{2k-1}\right)\,dx+
\left(B_{k-1}+xH_{2k-1}\right)\,dy.
\]
By \cite[Corollary~3.4]{lizarbe2023density}, $\F$ is a Riccati foliation and
admits an invariant algebraic curve of degree one.

Assume now that $p\neq[0:0:1]$. Since the automorphism group of
$\Pj(1,1,k)$ acts transitively on the smooth locus
$\Pj(1,1,k)\setminus\{[0:0:1]\}$, we may assume that $p=[1:0:0]$.
Let $\omega=A\,dx+B\,dy+C\,dz$ be a quasi-homogeneous $1$-form defining $\F$, where
\[
\deg(A)=2k,\qquad
\deg(B)=3k,\qquad
\deg(C)=2k+1.
\]
Since $p=[1:0:0]$ is the unique singularity of $\F$, we may write
\[
\begin{aligned}
B&=B_0x^{2k}+B_1x^{2k-1}+\cdots+B_{2k},\\
C&=C_0x^{k+1}+C_1x^{k}+\cdots+C_{k+1},
\end{aligned}
\]
where each $B_i$ and $C_i$ is a quasi-homogeneous polynomial in the variables
$y$ and $z$.
The Euler condition $xA+yB+kzC=0$ implies $yB_{2k}+kzC_{k+1}=0$.
Hence there exists a homogeneous polynomial $H_k$ of degree $k$ such that $B_{2k}=kzH_k$ and $C_{k+1}=-yH_k$.
Since $x$ is the unique singularity,
Proposition~\ref{P:Q-milnor} gives $\mu^\Q(\F,p)=\mu^\Q(\F)=2k+1$.
Passing to the affine chart centered at $p$, the induced $1$-form becomes $\omega_0=\widetilde B\,dy+\widetilde C\,dz$, where $\widetilde B=B_1+B_2+\cdots+B_{2k-1}+kzH_k$ and $\widetilde C=C_1+C_2+\cdots+C_k-yH_k$.
Expanding the homogeneous parts, we obtain
\[
\begin{aligned}
\widetilde B
={}&
b_1y+b_2y^2+\cdots+b_ky^k
+a_1z+a_2yz+\cdots+a_ky^{k-1}z
\\
&
+b_{k+1}y^{k+1}
+\cdots
+b_{2k-1}y^{2k-1}
+kz(my^k+nz),
\\
\widetilde C
={}&
c_1y+c_2y^2+\cdots+c_ky^k
+ez-y(my^k+nz),
\end{aligned}
\]
where $m,n,e,a_i,b_i,c_i\in\C$.
Therefore,
\begin{equation}\label{eq:local-form}
\omega_0=
\widetilde B\,dy+\widetilde C\,dz.
\tag{1}
\end{equation}    

\smallskip

\noindent\textbf{Case 1.} Assume that $e=0$.

If $n=0$, then $\widetilde C=-y^mh(y)$, where $1\le m\le k$ and $h(0)\neq0$. Hence
\[
\mu(\F,p)
=
\dim_\C
\frac{\C\{y,z\}}
{(\widetilde B,y^mh(y))}
=
\dim_\C
\frac{\C\{y,z\}}
{(\widetilde B,y^m)}
\le k,
\]
which contradicts $\mu(\F,p)=2k+1$. Therefore $n\neq0$. After multiplying $\omega$ by a nonzero constant, we may assume $n=1$.
Now write
\[
\widetilde C
=
-y\left(
-\sum_{i=1}^{k}c_iy^{i-1}
+my^k
+z
\right).
\]
Performing the change of coordinates
\[
z\longmapsto
z-\sum_{i=1}^{k}c_iy^i+my^{k+1},
\]
we obtain $\widetilde C=-yz$.
Substituting this expression into \eqref{eq:local-form}, we obtain
\[
\omega_0
=
\left(
\sum_{i=1}^{2k-1}\widetilde b_i y^i
+
z\sum_{j=1}^{k}\widetilde a_jy^{j-1}
+kz^2
\right)dy
-yz\,dz,
\]
where
\[
\widetilde a_j=
\begin{cases}
a_1+2kc_1,& j=1,\\
a_j+2kc_j-(j-1)c_{j-1},&2\le j\le k.
\end{cases}
\]
Since $\mu(\F,p)=2k+1$, it follows that
\[
\widetilde a_j=0,\qquad
\widetilde b_i=0,\quad
1\le i\le 2k-2,
\]
and $\widetilde b_{2k-1}\neq0$.
Finally, after the rescaling
\[
y\longmapsto\lambda y,
\qquad
z\longmapsto\lambda^{2k-2}z,
\]
for a suitable $\lambda\in\C^*$, we may normalize $\widetilde b_{2k-1}=1$. Therefore,
\[
\omega_0=
\left(
y^{2k-1}
+
z\sum_{j=2}^{k}a_jy^{j-1}
+kz^2
\right)dy
-yz\,dz,
\]
which is the normal form in {\rm(2a)}.

\smallskip

\noindent\textbf{Case 2.} Assume that $e\neq0$.

After multiplying $\omega$ by a nonzero constant, we may assume $e=1$.
Replacing
\[
z\longmapsto
z+\sum_{i=1}^{k}\alpha_i y^i,
\]
for suitable constants $\alpha_i$, we may eliminate the terms
$c_1,\ldots,c_k$ in $\widetilde C$. Thus $\widetilde C=z-my^{k+1}-nyz$, where $m,n\in\C$.

We distinguish two cases.

\medskip

\noindent
\textbf{Case 2.1.} $n=0$.

If $m=0$, then $\widetilde C=z$, and therefore
\[
\mu(\F,p)
=
\dim_\C
\frac{\C\{y,z\}}
{(\widetilde B,z)}
=
\dim_\C
\frac{\C\{y\}}
{(\widetilde B(y,0))}
\le 2k,
\]
contradicting $\mu(\F,p)=2k+1$. Hence $m\neq0$, and after rescaling $z$ we may assume $m=1$. Thus $\widetilde C=z-y^{k+1}$.
Since $\mu(\F,p)=2k+1$, it follows that $b_i=0$ for $1\le i\le k$, and $a_i+b_{k+i}=0$ for $1\le i\le k-1$.
Therefore,
\[
\omega_0=
\left(
(z-y^{k+1})
\sum_{j=1}^{k-1}a_jy^{j-1}
+kzy^k
\right)dy
+
(z-y^{k+1})\,dz,
\]
which is the normal form {\rm(2b)}.

\medskip

\noindent
\textbf{Case 2.2.} $n\neq0$.

After rescaling $z$, we may assume $n=1$. Hence $\widetilde C=z-my^{k+1}-yz$.
If $m=0$, then $\widetilde C=z(1-y)$, and
\[
\mu(\F,p)
=
\dim_\C
\frac{\C\{y,z\}}
{(\widetilde B,z(1-y))}
=
\dim_\C
\frac{\C\{y\}}
{(\widetilde B(y,0))}
\le 2k,
\]
contradicting $\mu(\F,p)=2k+1$. Therefore $m\neq0$, and after rescaling $z$ we may assume $m=1$. Thus $\widetilde C=z-y^{k+1}-yz$.
Since
\[
z=\frac{y^{k+1}}{1-y}
=
y^{k+1}(1+y+y^2+\cdots),
\]
the equality $\mu(\F,p)=2k+1$ implies that the coefficients of the pull-back of $\widetilde B$ under the
parametrization
\[
y\longmapsto
\left(y,\frac{y^{k+1}}{1-y}\right)
\]
vanish up to order $2k$. Equivalently, $b_i=0$ for $1\le i\le k$, and $b_{k+j}+a_j+\cdots+a_1=0$ for $1\le j\le k-1$.
Substituting these relations into $\widetilde B$ gives
\[
\widetilde B
=
\sum_{j=1}^{k-1}
a_jy^{j-1}(z-y^{k+1}-yz)
+kz(y^k+z).
\]
Therefore,
\[
\omega_0=
\left(
(z-y^{k+1}-yz)
\sum_{j=1}^{k-1}a_jy^{j-1}
+kz(y^k+z)
\right)dy
+
(z-y^{k+1}-yz)\,dz,
\]
which is the normal form {\rm(2c)}.
This completes the proof.
\end{proof}

Finally, we classify the foliations of normal degree $2k+2$ on
$\Pj(1,1,k)$, where $k>1$, having a unique singularity.

\begin{thm}\label{T:case2k+2(1,1,k)}
Let $\F$ be a foliation on $\Pj(1,1,k)$ of normal degree $2k+2$, with a unique singularity at $p=[0:0:1]$.
Then exactly one of the following holds:
\begin{enumerate}
    \item $m_p(\F)=k+1$, and $\F$ is a Riccati foliation. In particular, $\F$ admits at least one invariant curve; or

    \item $m_p(\F)=1$, and $p$ is a saddle-node singularity. In this case, $\F$ admits no invariant curve.
\end{enumerate}
\end{thm}

\begin{proof}
Since
\[
m_p(\mathcal F)\equiv \deg(N_{\mathcal F})-1
\equiv 1 \pmod{k},
\]
we have
\[
m_p(\mathcal F)=1
\quad\text{or}\quad
m_p(\mathcal F)=k+1.
\]

Assume first that $m_p(\mathcal F)=k+1$. Then, in suitable local
coordinates around $p=[0:0:1]$, the foliation is induced by
\[
\omega_0=
\left(A_{k+1}-yH_{2k}\right)dx+
\left(B_{k+1}+xH_{2k}\right)dy.
\]
By \cite[Corollary~3.4]{lizarbe2023density}, $\F$ is a Riccati foliation and
admits an invariant algebraic curve of degree one.

Assume now that $m_p(\mathcal F)=1$. Then $p$ is either a saddle-node
or a nilpotent singularity. We show that the latter is impossible.
Indeed, if $p$ is nilpotent, then
\[
\omega_0=
\left(A_{k+1}-yH_{2k}\right)dx+
\left(y+B_{k+1}+xH_{2k}\right)dy,
\]
where
\[
A_{k+1}
=a_0x^{k+1}+a_1x^ky+\cdots+a_{k+1}y^{k+1}.
\]
Let
\[
P=A_{k+1}-yH_{2k},
\qquad
Q=y+B_{k+1}+xH_{2k}.
\]
Since the initial form of $Q$ is $y$, whereas the initial form of $P$ is
$A_{k+1}$, we have $(P,Q)_{(0,0)}=\nu_{(0,0)}(P)\nu_{(0,0)}(Q)$ if and only if the initial forms of $P$ and $Q$ have no common factor.
If $a_0\neq0$, then
\[
(P,Q)_{(0,0)}
=1\cdot(k+1)
=k+1
<k\mu^\Q(\F,p)=(k+1)(2k+1)+k,
\]
a contradiction. Therefore $a_0=0$.
Hence $A_{k+1}=yA_k$, and consequently $P=yP_1$. It follows that the curve $\{y=0\}$ is invariant under $\mathcal F$, contradicting Theorem~\ref{T:(1,1,k)}.
Therefore the nilpotent case cannot occur. Hence $p$ is a saddle-node singularity. The final assertion follows from Theorem~\ref{T:(1,1,k)}.
\end{proof}
We now turn to the second case of Theorem~6.4
when $k=2$, and give an explicit
classification, up to automorphisms of
$\mathbb{P}(1,1,2)$, of foliations $\mathcal{F}$ of normal
degree $6$ having a unique singularity,
which $m_p(\mathcal{F})=1$ and $p=[0:0:1]$  is a saddle-node.
In suitable local coordinates $(x,y)$ centered
at $p$, the foliation $\mathcal{F}$ is represented by
the $1$-form
\[
\omega=(A_3+yG)\,dx+(x+B_3-xG)\,dy,
\]
where
 $A_3=a_{30}x^3+a_{21}x^2y+a_{12}xy^2+a_{03}y^3$ , $B_3=b_{30}x^3+b_{21}x^2y+b_{12}xy^2+b_{03}y^3$ and  $G$ is an homogeneous polynomial of degree 4. Since, by Theorem~6.4, $\mathcal{F}$ admits
no invariant algebraic curve, neither of
the coordinate axes is invariant. Hence,
\[
a_{30}b_{03}\neq 0.
\]  
Moreover, there exists an automorphism of $\mathbb P^2$
preserving the coordinate axes $x=0$ and $y=0$
and mapping the tangent line $L_{p_0}$ at a generic
point $p_0=(0,y_0)$ of a leaf of $\mathcal F$
to the line at infinity $z=0$, while sending $p_0$
to $[0:1:0]$.

Consequently, up to such an automorphism, we may
assume that $G=xG_3$, where $G_3=g_{30}x^3+g_{21}x^2y+g_{12}xy^2+g_{03}y^3$.
Thus, $\mathcal{F}$ is represented by
\begin{equation}\label{eq:omega-saddle}
	\omega=(A_3+xyG_3)\,dx
	+(x+B_3-x^2G_3)\,dy.
\end{equation}

Writing
\[
A=A_3+xyG_3,
\qquad
B=x+B_3-x^2G_3,
\]
the condition that $p$ is the unique
singularity of $\mathcal{F}$ is equivalent to
\[
I_0(A,B)=17,
\]
where $I_0(A,B)$ denotes the local
intersection multiplicity at the origin.
We obtain the following classification.

\begin{prop}\label{prop:saddle-node}
	Every foliation of normal degree six on $\mathbb P(1,1,2)$
having a unique singularity of saddle-node type is,
up to an automorphism of $\mathbb P(1,1,2)$, precisely
a member of the family $\{\mathcal F_\lambda\}_{\lambda\in\mathbb C}$
defined by the $1$-forms
\begin{equation}\label{eq:normal-form-sn}
\omega_\lambda
=(x^3+\lambda x^4y+x^2y^3)\,dx
+(x+\lambda x^2y+y^3-\lambda x^5-x^3y^2)\,dy.
\end{equation}
\end{prop}
\begin{proof}
A straightforward computation shows that
\[
I_{0}(A,B)=17,
\]
where $Ax^3+\lambda x^4y+x^2y^3 $ and $ x+\lambda x^2y+y^3-\lambda x^5-x^3y^2$.
Therefore, it only remains to prove the necessity. Assume that the foliation $\mathcal{F}$ is given by~\eqref{eq:omega-saddle}
\[\omega=(A_3+xyG_3)\,dx
	+(x+B_3-x^2G_3)\,dy,\]
    and $I_0(A,B)=17$.
Since the linear part of $B$ is $x$, the curve $B=0$ is smooth.
Moreover, because $b_{03}\ne 0$, it admits a parametrization
\[
y=t,
\qquad
x=x(t)=-b_{03}t^3+O(t^4),
\]
Consequently, $I_0(A,B)=3$, if $a_{03}\ne 0$, a contradiction. Hence $a_{03}=0$. Similarly, if $a_{12}\ne 0$ then $I_0(A,B)=5$. Therefore $a_{12}=0$, and
\[
A_3=a x^3+c x^2y,
\]
where
\[
a=a_{30}\ne 0,
\qquad
c=a_{21}.
\]

We now use an automorphism involving the coordinate of weight $2$.
Let $[X:Y:Z]$ be weighted homogeneous coordinates on
$\mathbb P(1,1,2)$. For suitavel constant $r,a$ consider the quadratic polynomial
$
q(X,Y)=rX^2+sY^2,
$
and the automorphism
\[
[X:Y:Z]\longmapsto [X:Y:Z-q(X,Y)].
\]

On the local chart centered at $[0:0:1]$, it can
be written as
\[
\Phi(u,v)=\left(\frac{u}{\sqrt{ 1-q(u,v)}},\frac{v}{\sqrt{1-q(u,v)}}\right)=\left(\frac{u}{\sqrt h},\frac{v}{\sqrt h}\right),
\]
where the square root is taken locally with value $1$ at the origin.
Since \(A_3(u,v)=u^2(au+cv)\), a direct computation yields

\[
\begin{aligned}
h^3\Phi^*\omega
={}&\bigl(\widetilde A_3+uv\widetilde G_3\bigr)\,du\\
&+\bigl(u+\widetilde B_3-u^2\widetilde G_3\bigr)\,dv,
\end{aligned}
\]
where
\[
\begin{aligned}
\widetilde A_3&=A_3+ru^2v,\\
\widetilde B_3&=B_3-2ru^3-suv^2,\\
\widetilde G_3&=G_3-suv(au+cv)+rB_3-ruq.
\end{aligned}
\]
In particular, \(\widetilde G_3\) is homogeneous of
degree three, so the transformed form has the
same structure as \(\omega\). Moreover, its cubic
coefficients satisfy
\[
\begin{gathered}
a'_{30}=a_{30},\qquad a'_{21}=a_{21}+r,\\
b'_{30}=b_{30}-2r,\qquad b'_{21}=b_{21},\\
b'_{12}=b_{12}-s,\qquad b'_{03}=b_{03}.
\end{gathered}
\]
In particular, choosing
\[
r=-a_{21},
\qquad
s=b_{12},
\]
we obtain
\[
a_{21}'=b_{12}'=0.
\]

Renaming the coordinates and coefficients, we may therefore assume
that
\[
A=a x^3+xyG_3,
\qquad
B=x+B_3-x^2G_3,
\]
where
\[
B_3=d x^3+e x^2y+b y^3,
\qquad
ab\ne 0.
\]

Therefore the condition $I_0(A,B)=17$  is equivalent to the condition
\begin{equation}\label{eq:simply1}
	I_0(ax^2+yG_3,B)=14=I_0(ax^2+yG_3-axB,B).
\end{equation}

Since, the curve $B=0$ admits a parametrization \[
y=t,
\qquad
x=x(t)=-bt^3+O(t^4),
\] 
~\eqref{eq:simply1} is equivalent to the multiplicity of
\begin{equation}\label{eq:simply2}
	R(x(t),t)= d_2 t^2+d_3t^3+\cdots+d_{14}t^{14}+\cdots,
\end{equation}
where $R=ax^2+yG_3-axB$ is 14. So we must find conditions on the coefficients of $A$ and $ B$ such that
$d_j=0$ for all $j=2, \cdots,13$ and $d_{14} \neq 0$
Write
\[
G_3=g_{30}x^3+g_{21}x^2y+g_{12}xy^2+g_{03}y^3.
\]
Then
\[
R=g_{03}y^4
 +(g_{12}-ab)xy^3
 +g_{21}x^2y^2
+(g_{30}-ae)x^3y
 -ad x^4
 +ax^3G_3
 \]
and the coefficient of $R(x(t),t)$ the orders $2,3,5,7,9,11,13$ of are zero and the orders $4,6,8,10$ and $12$ yields
\[
\begin{aligned}
d_4=0 &\implies g_{03}=0,\\
d_6=0 &\implies g_{12}=ab,\\
d_8=0 &\implies g_{21}=0,\\
d_{10}=0 &\implies g_{30}=ae,\\
d_{12}=0 &\implies d=0.
\end{aligned}
\]
Therfore we have
\[
G_3=aex^3+abxy^2,
\qquad
R=ax^3G_3.
\]
Consequently,
\[
R(x(t),t)=a^2b^5t^{14}+O(t^{18}).
\]
Since \(ab\neq 0\), it follows that
\(I_0(R,B)=14\), as required.
Hence the form is
\[
\begin{aligned}
\omega={}&
\left(ax^3+ae x^4y+ab x^2y^3\right)\,dx\\
&+
\left(x+e x^2y+b y^3
      -ae x^5-ab x^3y^2\right)\,dy.
\end{aligned}
\]
It remains to normalize the nonzero constants $a$ and $b$.
Consider the diagonal change of coordinates
\[
x=\mu X,\qquad y=\nu Y,\qquad \mu\nu\neq 0.
\]
After dividing the resulting $1$-form by $\mu\nu$, the
coefficients of $X^3dX$ and $Y^3dY$ become
\[
\frac{a\mu^3}{\nu}
\quad\text{and}\quad
\frac{b\nu^3}{\mu},
\]
respectively. Choosing
\[
\nu=a\mu^3,\qquad
\mu^8=\frac{1}{a^3b},
\]
both coefficients are equal to $1$. Moreover,
$ab\mu^2\nu^2=1$. Setting $\lambda=e\mu\nu$,
we obtain the resulting form is precisely \eqref{eq:normal-form-sn}.
\end{proof}

The next three examples exhibit foliations without invariant algebraic
curves, obtained by applying Theorem~\ref{T:(1,1,k)}.

\begin{example}
Consider the foliation on $\Pj(1,1,k)$ defined by
\[
\begin{aligned}
\omega
={}&
\left(
-kz^2-x^{k-1}yz+x^{k-2}y^{k+2}
\right)dx
\\
&
+\left(
x^kz-x^{k-1}y^{k+1}+kxzy^k
\right)dy
\\
&
-(xz-y^{k+1})\,dz.
\end{aligned}
\]
It has $\deg(N_{\mathcal F})=2k+1$.
Moreover,
\[
\operatorname{Sing}(\mathcal F)=\{p\},
\qquad
p=[1:0:0],
\]
where $p$ is a saddle-node singularity.
In the affine chart centered at $p$, the induced $1$-form is
\[
\omega_0=
(z-y^{k+1}+kzy^k)\,dy
+
(z-y^{k+1})\,dz.
\]
A direct computation yields $\mu^\Q(\mathcal F,p)=2k+1$. By Proposition~\ref{P:Q-milnor}, $\mu^\Q(\mathcal F)=2k+1$. Hence $\mu^\Q(\mathcal F,p)=\mu^\Q(\mathcal F)$, and therefore $p$ is the unique singularity of $\mathcal F$, and it is a saddle-node.
\end{example}

\begin{example}
Consider the affine chart $z=1$ of $\Pj(1,1,k)$. The foliation is defined by
\[
\omega=
x^k\bigl(-x+y^{k+1}\bigr)\,dx
+
\bigl(x+y^{k+1}-x^{k+1}y^k\bigr)\,dy.
\]
This affine $1$-form induces a foliation $\mathcal F$ on $\Pj(1,1,k)$ with $\deg(N_{\mathcal F})=2k+2$.
Moreover,
\[
\operatorname{Sing}(\mathcal F)=\{p\},
\qquad
p=[0:0:1],
\]
where $p$ is a saddle-node singularity.
Indeed, a direct computation in the affine chart $z=1$ yields
\[
\mu^{\Q}(\mathcal F,p)
=
\frac{(2k+1)(k+1)+k}{k}.
\]
According to Proposition~\ref{P:Q-milnor},
\[
\mu^{\Q}(\mathcal F)
=
\frac{(2k+1)(k+1)+k}{k}.
\]
Hence $\mu^{\Q}(\mathcal F,p)=\mu^{\Q}(\mathcal F)$, and therefore $p$ is the unique singularity of $\mathcal F$, and it is a saddle-node.
\end{example}

\begin{example}
Consider the foliation on $\Pj(1,1,k)$ defined by
\[
\begin{aligned}
\omega
={}&
\left(
xz^3-y^{k+1}z^2-x^{3k}y
\right)dx
\\
&
-
\left(
xy^kz^2-y^{2k+1}z+x^{3k+1}
\right)dy
\\
&
-\frac1k
\left(
x^2z^2-2xy^{k+1}z-y^{2k+2}
\right)dz.
\end{aligned}
\]
This $1$-form defines a foliation $\mathcal F$ on $\Pj(1,1,k)$ with $\deg(N_{\mathcal F})=3k+2$.
Moreover,
\[
\operatorname{Sing}(\mathcal F)=\{p\},
\qquad
p=[0:0:1].
\]
Indeed, in the affine chart $z=1$, the induced $1$-form is
\[
\omega_0=
(x-y^{k+1}-x^{3k}y)\,dx
-
(xy^k-y^{2k+1}+x^{3k+1})\,dy.
\]
A direct computation yields
\[
\mu^{\Q}(\mathcal F,p)
=
\frac{6k(k+1)+1}{k}.
\]
On the other hand, Proposition~\ref{P:Q-milnor} gives
\[
\mu^{\Q}(\mathcal F)
=
\frac{6k(k+1)+1}{k}.
\]
Hence $p$ is the unique singularity of $\mathcal F$. Since the linear part of
$\omega_0$ at $p$ is $x\,dx$, the singularity $p$ is nilpotent.
\end{example}

\begin{lemma}\label{L:wellformed-c2}
Let
$a\le b\le c$ be pairwise coprime positive integers. Suppose that
$$ab+ac+bc>abc.$$ Then $a=1$ or $(a,b,c)=(2,3,5)$.
\end{lemma}

\begin{proof}
Dividing the inequality $ab+ac+bc>abc$ by $abc>0$, we obtain
\[
\frac1a+\frac1b+\frac1c>1.
\]

Suppose first that $a>1$. Since $1<a<b<c$, we have $a\ge2$, $b\ge3$ and $c\ge4$.
If $a\ge3$, then
\[
\frac1a+\frac1b+\frac1c
\le
\frac13+\frac14+\frac15
=
\frac{47}{60}<1,
\]
a contradiction. Hence $a=2$.
If $b\ge4$, then $c\ge5$, and
\[
\frac12+\frac1b+\frac1c
\le
\frac12+\frac14+\frac15
=
\frac{19}{20}<1,
\]
again a contradiction. Thus $b=3$.
Finally,
\[
\frac12+\frac13+\frac1c>1
\]
implies $\frac1c>\frac16$, and hence $c<6$. Since $c>3$, we have $c\in\{4,5\}$. The condition $\gcd(a,c)=1$ excludes $c=4$, since $\gcd(2,4)=2$. Therefore $c=5$, and hence $(a,b,c)=(2,3,5)$.
\end{proof}

\begin{cor}\label{C:P1bc}
Let
$
X=\mathbb P(1,b,c)$, where
$ 1<b<c$ and 
$\gcd(b,c)=1$, 
and let $\F$ be a foliation on $X$ with a unique singularity $x$. Write $N_{\F}\equiv dH$, where $H=\mathcal O_{\mathbb P(1,b,c)}(1)$ and
$d=\deg(N_{\F})$. Assume that the algebraic multiplicity of $\F$ at $x$
is one.
Then $\F$ admits no invariant algebraic curve under the following
conditions:

\[
d>
\frac23
\left(
1+b+c+
\sqrt{(b+c)^2-(b+c)+1}
\right)
\]
and
\[
d\ge
1+b+c+\sqrt{(c-b)^2+1},
\qquad \hbox{if $x$ is smooth}.
\]

If $x=[0:1:0]$, then
\[
d\ge
1+b+c+\sqrt{b^2+(c-1)^2}.
\]

If $x=[0:0:1]$, then
\[
d\ge
1+b+c+\sqrt{c^2+(b-1)^2}.
\]
\end{cor}

\begin{cor}\label{C:P235}
Let $X=\mathbb P(2,3,5)$, and let $\F$ be a foliation on $X$ with a unique singularity $x$.
Write $N_{\F}\equiv dH$, where $H=\mathcal O_{\mathbb P(2,3,5)}(1)$
and $d=\deg(N_{\F})$. Assume that the algebraic multiplicity of $\F$ at
$x$ is one. If $d\ge 17$, then $\F$ admits no invariant algebraic curve. 
\end{cor}

\begin{proof}
    This follows directly from Theorem \ref{T:classone} and Lemma
    \ref{L:particularlocalobstruction}.
\end{proof}

\bibliography{references}{}
\bibliographystyle{plain}

\end{document}